\documentclass[11pt,amsfonts]{article}
\usepackage{graphicx}
\usepackage{latexsym}
\usepackage{amssymb}
\usepackage{amsmath}
\usepackage[latin1]{inputenc}
\usepackage{layout}
\newtheorem{prop}{Proposition}
\newtheorem{lemma}{Lemma}

\newtheorem{theorem}{Theorem}
\newtheorem{remark}{Remark}

\def\real{{\mathord{{\rm I\kern-2.8pt R}}}}        
\def\inte{{\mathord{{\rm I\kern-2.8pt N}}}}

\def\sZZ{{\rm Z\kern-2.8ptem{}Z}}

\def\z{{\mathchoice
  {\sZZ}
  {\sZZ}
  {\rm Z\kern-0.30em{}Z}
  {\rm Z\kern-0.25em{}Z} }}
\def\sQQ{{\kern 0.27em \vrule height1.45ex width0.03em depth0em
          \kern-0.30em \rm Q}}
\def\qu{{\mathchoice
    {\sQQ}
    {\sQQ}
  {\kern 0.225em \vrule height1.05ex width0.025em depth0em \kern-0.25em \rm Q}
  {\kern 0.180em \vrule height0.78ex width0.020em depth0em \kern-0.20em \rm Q}
        }}
\def\sCC{{\kern 0.27em \vrule height1.45ex width0.03em depth0em
          \kern-0.30em \rm C}}
\def\complex{{\mathchoice
    {\sCC}
    {\sCC}
  {\kern 0.225em \vrule height1.05ex width0.025em depth0em \kern-0.25em \rm C}
  {\kern 0.180em \vrule height0.78ex width0.020em depth0em \kern-0.20em \rm C}
        }}

\newcommand{\ba}{\begin{array}}
\newcommand{\ea}{\end{array}}
\newcommand{\be}{\begin{equation}}
\newcommand{\ee}{\end{equation}}
\newcommand{\bea}{\begin{eqnarray}}
\newcommand{\eea}{\end{eqnarray}}
\newcommand{\beaa}{\begin{eqnarray*}}
\newcommand{\eeaa}{\end{eqnarray*}}

\newcommand{\eps}{\varepsilon}

\def\z{\zeta}

\font\tenmath=msbm10 \font\sevenmath=msbm7 \font\fivemath=msbm5
\newfam\mathfam \textfont\mathfam=\tenmath
\scriptfont\mathfam=\sevenmath \scriptscriptfont\mathfam=\fivemath

\def \={{\buildrel {\rm (law)} \over =}}

\def\cd·
\def\cds{\cdots}

\def\qed{ \hfill \vrule width.25cm height.25cm depth0cm\smallskip}

\newcommand{\basa}{\begin{assumption}}
\newcommand{\easa}{\end{assumption}}

\newcommand{\bas}{\begin{assum}}
\newcommand{\eas}{\end{assum}}

 \def\cd·
\def\cds{\cdots}

\newcommand{\ignore}[1]{}
\begin{document}

\renewcommand{\thefootnote}{\fnsymbol{footnote}}

\date{ }
\title{ Approximation of the finite dimensional distributions of  multiple fractional integrals }
\author{Xavier Bardina$^{1}, $\thanks{The author is supported by the grant MEC-FEDER Ref. MTM2006-06427.}\qquad  Khalifa Es-Sebaiy $^{2}\qquad $%
Ciprian A. Tudor $^{3}$, \footnote{Associate member of the team Samos, Universit\'e de Panth\'eon-Sorbonne Paris 1 }\vspace*{0.1in} \\
$^{1}$Departament de Matem\`atiques, Universitat Aut\`onoma de Barcelona\\
08193-Bellaterra (Barcelona), Spain. \\
bardina@mat.uab.cat \vspace*{0.1in} \\
$^{2}$SAMOS-MATISSE, Centre d'Economie de La Sorbonne, \\
Université de Paris 1 Panth\'eon-Sorbonne,\\
90, rue de Tolbiac, 75634, Paris, France. \\
Khalifa.Es-Sebaiy@malix.univ-paris1.fr \vspace*{0.1in} \\
 $^{3}$ Laboratoire Paul Painlev\'e, Universit\'e de Lille 1\\
 F-59655 Villeneuve d'Ascq, France.\\
 \quad tudor@math.univ-lille1.fr\vspace*{0.1in}}
\maketitle

\begin{abstract}
We construct a family $I_{n_{\eps}}(f)_{t}$ of continuous stochastic processes that converges in the sense of finite dimensional distributions
 to a multiple Wiener-It\^o  integral $I_{n}^{H}(f1^{\otimes n}_{[0,t] })$ with respect to the fractional Brownian motion. We assume that $H>\frac{1}{2}$ and we prove our approximation result for the integrands $f$ in a rather general class.
\end{abstract}

\vskip0.5cm

{\bf  2000 AMS Classification Numbers:} 60B10, 60F05, 60H05.

 \vskip0.3cm


\section{Introduction}
A first result concerning the approximation of iterated stochastic
integrals has been given in \cite{Avram}. Consider
$(X^{\eps})_{\eps>0}$  a family of semimartingales with paths in
the Skorohod space  $\mathcal D([0,1])$ that converges weakly in
this space to another semimartingale $X$, as $\eps$ tends to zero.
It has been proven in \cite{Avram} that the couple
$(X^\varepsilon,[X^\varepsilon,X^\varepsilon])$ converges weakly in
${\cal{D}}([0,1]) $ as $\eps \to 0$ to  the couple  $(X,[X,X])$
($[X,X]$ denotes the usual semimartingale bracket) if and only if
for every $m\geq 1$ the vector $
(J_1(X^\varepsilon),\dots,J_m(X^\varepsilon))$ converges weakly in
${\cal{D}}([0,1]) $ as $\eps \to 0$ to  the vector
$(J_1(X),\dots,J_m(X))$. Here $J_{1} (X)_{t} = X_{t} $ and for
$k\geq 2, $ $J_{k}(X)_{t}=\int_{0}^{t} J_{k-1} (X)_{s-}dX_{s}$ (and
similarly for $J_k(X^\varepsilon)$). This result shows that in order
to obtain (joint) weak convergence of iterated Itô integrals we need
the convergence of $X^{\eps}$ to $X$ but also the convergence of the
second order variations. When our semimartingale is the Wiener
process, there are many  examples of families of processes with
absolutely continuous paths converging weakly to it in the topology
of $\mathcal C([0,1])$. In this case it is obvious that we do not have
convergence of the quadratic variations of such families to the
quadratic variations of the Brownian motion. This led to the problem of approximating iterated stochastic integrals with respect to the Browniam motion and later, with respect to  the fractional Brownian motion.

Let us recall some relatively recent results concerning the
approximation of iterated integrals with respect to a standard
Brownian motion by a family of processes with  continuous paths.
Consider a family of stochastic processes  $(\rho_{\eps})_{\eps >0}$
of the form
\begin{equation*}
\rho_{\eps} (t)= \int_{0}^{t} \theta _{\eps} (s) ds
\end{equation*}
such that $(\rho_{\eps}) _{\eps >0}$ converges weakly in $C_{0}([0,1])$ (the space of continuous function on $[0,1]$ which are null at zero) to the Wiener process. We will discuss two main examples: the case when $\theta_{\varepsilon}(s):=\frac{1}{\varepsilon}\sum_{k=1}^{\infty}
 \xi_k I_{[k-1,k)}\left(\frac{s}{\varepsilon^2}\right),$
 where $\{\xi_k\}$ is a sequence of independent, identically distributed
 random variables satisfying $E(\xi_1)=0$ and $\mbox{Var}(\xi_1)=1$ (these kernels are traditionally called {\it Donsker kernels}) and the case when $$\theta_{\varepsilon}(x):=\frac1{\varepsilon}(-1)^{N\left(\frac{x}{\varepsilon^2}\right)},$$
where $N=\{N(s);\,s\geq0\}$ is a standard Poisson process (these
kernels are usually called {\it Stroock kernels} or {\it Kac-Stroock
kernels} because they were introduced by Kac in \cite{kac} and used
by Stroock, \cite{stroock}, in order to obtain weak approximations
of the Brownian motion).
 In \cite{BJ}  the authors proved that, for a suitable function $f$ defined on $[0,1]^{\otimes n}$,  the family of multiple integrals $(I^{1}_{n_{\eps}}(f))_{\eps >0}$ with respect to $\rho  _{\eps} $ given by
 \begin{eqnarray}
  I^{1}_{n_{\eps}}(f)_{t}&=& \int_{[0,t] ^{ n} }f(t_{1},\dots,t_{n})d\rho_{\eps}(t_{1})\cdots d\rho_{\eps}(t_{n})\nonumber\\
  &=&\int_{[0,t] ^{ n} }f(t_{1},\dots,t_{n})\theta_{\eps }(t_{1})\cdots\theta_{\eps }(t_{n}) dt_{1}\cdots dt_{n}\label{i11}
  \end{eqnarray}
  converges weakly in $C_{0}([0,1])$ to the $n$th multiple Stratonovich integral of $f1_{[0,t]}^{\otimes n}$ with respect to the standard Brownian motion.
This is somehow expected because the Stratonovich integral usually
satisfies the differential rules  of the deterministic calculus. In
order to obtain as a limit a multiple It\^o integral (which has zero
expectation) one needs to subtract the ``trace" of
$I^{1}_{n_{\eps}}(f)$, that means, to suppress the values on the
diagonals. This following result has been obtained  in \cite{BJT2}:
for any $f\in L^{2}([0,1] ^{n})$, the family
$(I^{2}_{n_{\eps}}(f))_{\eps >0}$ given by
  \begin{eqnarray}\label{i2}
I^{2}_{n_{\eps}}(f)_{t}&=&\int_{[0,t]^n}f(x_1,x_2,\dots,x_n)\prod^n_{\footnotesize{\begin{array}{c}i,j=1\\i\neq
j\end{array}}}
I_{\{|x_i-x_j|>\varepsilon\}}d\rho_{\varepsilon}(x_1)\cdots\cdots
d\rho_{\varepsilon}( x_n)\nonumber\\
&=& \int_{[0,t]^n}f(x_1,x_2,\dots,x_n)\prod_{i=1}^n\theta_{\varepsilon}(x_i)\prod^n_{\footnotesize{\begin{array}{c}i,j=1\\i\neq
j\end{array}}}I_{\{|x_i-x_j|>\varepsilon\}}dx_1\cdots
dx_n,\end{eqnarray}
  converges weakly, in the sense of finite dimensional distributions (and in $C_{0}([0,1])$ for $n=2$), to the $n$th multiple It\^o integral $I_{n}(f1_{[0,t]}^{\otimes n})$.
Let us consider now the problem of approximating the fractional
Brownian motion $(B^{H}_{t})_{t\in [0,1]}$ and the multiple
integrals with respect to it. Recall that  the fractional Brownian
motion is a centered Gaussian process with covariance $R(t,s)=
\frac{1}{2}(t^{2H}+ s^{2H}-\vert t-s\vert ^{2H} )$ with $H\in
(0,1)$. It can be also expressed as a Wiener integral with respect
to a Wiener process $W$ by $B^{H}=\int_{0}^{t} K_{H}(t,s)dW_{s}$
where $K_{H}$ is a deterministic kernel defined on the set
$\{0<s<t\}$ and given by
\begin{equation}\label{K}
K_H(t,s)=c_H (t-s)^{H-\frac12}+ c_H\left(\frac12-H\right)\int_s^t
(u-s)^{H-\frac32}\left(1-\left(\frac{s}{u}\right)^{\frac12-H}\right)
du,
\end{equation} where $c_H$ is the  normalizing constant
 $ c_H=\left(\frac{2H\Gamma(\frac32-H)}{\Gamma(H+\frac12)
 \Gamma(2-2H)}\right)^{\frac12}.$  From this representation
and the weak  convergence of $\rho _{\eps} $ to $W$ it follows that
(see \cite{BJT1})  for any $H\in (0,1)$ the family of processes
$(\eta_{\eps})_{\eps >0}$ with
\begin{equation*}
\eta _{\eps } (t) = \int_{0}^{t} K_{H}(t,s) \theta _{\eps } (s)ds, \hskip0.5cm t\in [0,1]
\end{equation*}
converges weakly as $\eps \to 0$ in $C_{0}([0,1])$ to $B^{H}$. When $H>\frac{1}{2}$ the paths of $\eta _{\eps}$ are even absolutely continuous. Moreover, if $H>\frac{1}{2}$, the multiple integral with respect to $\eta_{\eps}$
\begin{equation}
\label{i3} I^{3}_{n_{\eps}}(f)_{t}= \int_{[0,t] ^{\otimes n}
}f(t_{1},\dots,t_{n})d\eta_{\eps}(t_{1})\cdots d\eta_{\eps}(t_{n})
\end{equation}
converges as $\eps \to 0$ in $C_{0}([0,1])$ to the multiple
Stratonovich integral of order $n$ of the function
$f1_{[0,t]}^{\otimes n}$ with respect to $B^{H}$. The purpose of
this work is to give an approximation result for the multiple
Wiener-It\^o integrals $I_{n}^{H}(f1_{[0,t]}^{\otimes n})$  with
respect to the fractional Brownian motion, for the integrand $f$ in
a rather general class of functions. Note that, as we recall in
Section 2, the multiple fractional integral $I_{n}^{H}$ can be
expressed as a multiple Wiener-It\^o integral with respect to the
Brownian motion. In fact, we have $I_{n}^{H}(f1_{[0,t]}^{\otimes
n})=I_{n}\left( \Gamma ^{(n)}_{H}f1_{[0,t]}^{\otimes n}\right)$
where $\Gamma ^{(n)}_{H}$ is a transfer operator. Concretely, we
show here that the family $(I_{n_{\eps}}(f))_{\eps >0}$ defined by
\begin{equation*}
I_{n _{\eps}}(f)_{t}= \int_{0}^{1}\ldots \int_{0}^{1} \left( \Gamma
^{(n)}_{H} f1_{[0,t]}^{\otimes n}\right)(x_{1}, \dots, x_{n})\left(
\prod_{i=1}^{n}\theta _{\eps }(x_{i}) \right) \prod_{i,j=1; i\not=
j}^{n} 1_{\{\vert x_{i}- x_{j} \vert >\eps \}}dx_{1}\ldots dx_{n}
\end{equation*}
converges, in the sense of finite dimensional distributions, to $(I_{n}^{H}(f1_{[0,t]}^{\otimes n}))_{t\in [0,1]}.$ Due to the rather complicate expression of the operator $\Gamma _{H}^{(n)}$ this result cannot be deduced from the result in \cite{BJT2} since the transfer principle for multiple fractional integrals actually implies that  $I_{n}^{H}(f1_{[0,t]}^{\otimes n})$ is equal to $I_{n}(g(t, \cdot )1_{[0,t] }^{\otimes n}) $ with some function $g$ depending on $f$. Because of the appearance  of the variable $t$ in the argument of $g$, the main result in \cite{BJT2} cannot be directly applied.  Another particularity of the multiple fractional integrals is that the expectation $EI_{1}^{H}(1_{A}) I_{1}^{H}(1_{B})$ is not zero when $A$ and $B$ are disjoint subsets of $[0,1]$ and this fact makes the proofs considerably more complex than in the standard Brownian motion case.

We structured our paper in the following way. Section 2 contains some preliminaries on multiple Wiener-It\^o integrals and multiple integrals with respect to the fractional Brownian motion. In Section 3 we prove our approximation result. We first regard the case when the integrand is a step function. We separated the case $n=2$ and $n\geq 3$ because in the first case the proof is less complex and more intuitive and it helps to understand the general case. Finally we extend our result from simple functions to a bigger class of functions.

\section{Preliminaries}
\subsection{Multiple Wiener-It\^o integrals}

In this paragraph we describe the basic elements of
calculus on Wiener chaos.  Let $(W_{t})_{t\in [0,1]}$ be a classical
Wiener process on a standard Wiener space $\left( \Omega
,{\mathcal{F}},\mathbf{P}\right) $. If $f\in
L^{2}([0,1]^{n})$ with $n\geq 1$ integer, we introduce the multiple Wiener-It%
ô integral of $f$ with respect to $W$. We refer to \cite{N} for
a detailed exposition of the construction and the properties of
multiple Wiener-Itô integrals.

Let $f\in {\mathcal{S}_{m}}$ be an elementary functions with $m$
variables  that can be written as \begin{equation*}
f=\sum_{i_{1},\ldots ,i_{m}}c_{i_{1},\ldots ,i_{m}}1_{A_{i_{1}}\times
\ldots \times A_{i_{m}}}
\end{equation*}%
where the coefficients satisfy $c_{i_{1},\ldots ,i_{m}}=0$ if two indices $%
i_{k}$ and $i_{l}$ are equal and the sets $A_{i}\in
{\mathcal{B}}([0,1])$ are pairwise-disjoints. For a such step
function $f$ we define

\begin{equation*}
I_{m}(f)=\sum_{i_{1},\ldots ,i_{m}}c_{i_{1},\ldots
,i_{m}}W(A_{i_{1}})\ldots W(A_{i_{m}})
\end{equation*}%
where we put $W([a,b])=W_{b}-W_{a}$. It can be seen that the application $%
I_{m}$ constructed above from ${\mathcal{S}_m}$ to $L^{2}(\Omega )$
is an
isometry on ${\mathcal{S}_m}$ , i.e.%
\begin{equation}
E\left[ I_{n}(f)I_{m}(g)\right] =n!\langle f,g\rangle
_{L^{2}([0,1]^{n})}\mbox{ if }m=n  \label{isom}
\end{equation}%
and%
\begin{equation*}
E\left[ I_{n}(f)I_{m}(g)\right] =0\mbox{ if }m\not=n.
\end{equation*}%
It also holds that
\begin{equation*}
I_{n}(f)=I_{n}\left( \tilde{f}\right)
\end{equation*}%
where $\tilde{f}$ denotes the symmetrization of $f$ defined by $\tilde{f}%
(x_{1},\ldots ,x_{n})=\frac{1}{n!}\sum_{\sigma \in S_{n}}f(x_{\sigma
(1)},\ldots ,x_{\sigma (n)})$.

Since the set ${\mathcal{S}_{n}}$ is dense in $L^{2}([0,1]^{n})$ for every $n\geq 1$ the mapping $%
I_{n}$ can be extended to an isometry from $L^{2}([0,1]^{n})$ to $%
L^{2}(\Omega)$ and the above properties hold true for this
extension. Note also that $I_{n}$ can be viewed as an iterated
stochastic integral
\begin{equation*}
I_{n}(f)=n!\int_{0}^{1}\int_{0}^{t_{n}}\ldots\int_{0}^{t_{2}}f(t_{1},%
\ldots,t_{n})dW_{t_{1}}\ldots dW_{t_{n}};
\end{equation*}
here the integrals are of Itô type; this formula is easy to show for
elementary functions $f$, and follows for general $f\in
L^{2}([0,1]^{n})$ by a density argument.

The product for two multiple integrals says that (see \cite{N}): if $f\in L^{2}([0,1]^{n})$ and
$g\in L^{2}([0,1]^{m})$ are symmetric functions, then
\begin{equation}
I_{n}(f)I_{m}(g)=\sum_{l=0}^{m\wedge
n}l!{{m}\choose{l}}{{n}\choose{l}}I_{m+n-2l}(f\otimes _{l}g)
\label{prod}
\end{equation}%
where the contraction $f\otimes _{l}g$ belongs to
$L^{2}([0,1]^{m+n-2l})$ for $l=0,1,\ldots ,m\wedge n$ and it is
given by
\begin{eqnarray}
&&(f\otimes _{l}g)(s_{1},\ldots ,s_{n-l},t_{1},\ldots ,t_{m-l})  \notag \\
&=&\int_{[0,1]^{l}}f(s_{1},\ldots ,s_{n-l},u_{1},\ldots
,u_{l})g(t_{1},\ldots ,t_{m-l},u_{1},\ldots ,u_{l})du_{1}\ldots
du_{l}. \label{contra}
\end{eqnarray}
When $l=0$, we will denote, throughout this paper,  by $f\otimes
g:=f\otimes _{0}g$.

\subsection{ Multiple fractional integrals}

Let us introduce here the multiple integrals with respect to the
fractional Brownian motion. We follow the approach in \cite{PAT}
(see also \cite{N} and \cite{Sa-To}). Let $f\in L^{1}([0,1] ^{n})$
and for every $0<\alpha <1$ define the operator

\begin{equation*}
\left( I_{t-}^{\alpha , n}f\right) (x_{1}, \ldots , x_{n})=
\frac{1}{\left(\Gamma (\alpha)\right) ^{n}  } \int_{x_{1}}^{t}
\ldots \int_{x_{n}}^{t} \frac{f(t_{1}, \ldots , t_{n})}{\prod
_{j=1}^{n} (t_j -x_{j}) ^{1-\alpha }}dt_{1}\ldots dt_{n}
\end{equation*}
for every $x_{1}, \dots, x_{n} \in [0,t]$ with $t\in[0,1]$.

We have the following properties:

\begin{description}
\item{$\bullet$} if $f= f_{1} \otimes \cdots \otimes f_{n} $ with $f_{i}\in L^{1}([0,1])$ then
\begin{equation*}
(I_{t-}^{\alpha , n}f ) (x_{1}, \ldots , x_{n}) = (I_{t-}^{\alpha , 1} f_{1})(x_{1} ) \ldots (I_{t-}^{\alpha , 1} f_{n})(x_{n} )
\end{equation*}
for every $x_{1}, \dots, x_{n} \in [0,t]$.
\item{$\bullet$} If $H>\frac{1}{2}$ then
\begin{equation}\label{i1}
c_{H} \Gamma \left( H+\frac{1}{2}\right) s^{\frac{1}{2} -H} \left(I_{1-} ^{H-\frac{1}{2}, 1} (x^{H-\frac{1}{2} }1_{[0,t] })\right)(s) =K_{H}(t,s)
\end{equation}
where $K_{H}$ is the standard kernel of the fractional Brownian motion (\ref{K}).
\end{description}

We introduce the space $\left|{\cal{H}}\right| ^{\otimes n}$  of measurable functions $f:[0,1] ^{n} \to \mathbb{R}$ such that
\begin{equation*}
\int_{[0,1] ^{2n}} \vert f(u_{1}, \ldots , u_{n}) f(v_{1}, \ldots ,
v_{n})\vert \left( \prod _{j=1}^{n} \psi (u_{j}, v_{j}) \right)
du_{1}\cdots du_{n}dv_{1}\cdots dv_{n}<\infty
\end{equation*}
where $\psi (s,t)= H(2H-1) \vert s-t\vert ^{2H-2}$.

\begin{remark}For any $H>\frac{1}{2}$ we have (see \cite{N}, \cite{PAT})
$$L^{2}([0,1] ^{n}) \subset L^{\frac{1}{H}}([0,1]^{n})\subset \left| {\cal{H}}\right| ^{\otimes n}.$$
\end{remark}

Define the operator $\Gamma ^{(n)}_{H} : \left| {\cal{H}}\right| ^{\otimes n} \to L^{2}([0,1]^{n})$
\begin{equation}
\label{opGamma} (\Gamma  ^{(n)}_{H}f)(t_{1}, \ldots , t_{n}) =\left[
c_{H} \Gamma \left( H+\frac{1}{2}\right)\right] ^{n} \prod
_{j=1}^{n} t_{j} ^{\frac{1}{2}-H}\left( I_{1-}^{H-\frac{1}{2},
n}\right) (f(x_{1}, \dots, x_{n})\prod_{j=1}^{n}x_{j}
^{H-\frac{1}{2}})(t_{1}, \ldots , t_{n})
\end{equation}
Then the operator $\Gamma ^{(n)}_{H}$ is an isometry between $\left| {\cal{H}}\right| ^{\otimes n} $ and  $ L^{2}([0,1]^{n})$ where we endow the space $\left| {\cal{H}}\right| ^{\otimes n} $ with the following inner product
\begin{equation*}
\langle f, g\rangle _{ {\cal{H}} ^{\otimes n} }=\int_{[0,1] ^{2n}}
f(u_{1}, \ldots , u_{n}) g(v_{1}, \ldots , v_{n}) \left( \prod
_{j=1}^{n} \psi (u_{j}, v_{j}) \right) du_{1}\cdots
du_{n}dv_{1}\cdots dv_{n}.
\end{equation*}

Note that
\begin{description}
\item{$\bullet$ } if $f=1_{[0,b]}$ then by (\ref{i1})
\begin{equation*}
(\Gamma ^{(1)}_{H} 1_{[0,b]} )(s)=c_{H} \Gamma \left( H+\frac{1}{2}\right)s^{\frac{1}{2}-H} I_{1-}^{H-\frac{1}{2}}(x^{H-\frac{1}{2}}1_{[0,b]})(s)=K_{H}(b,s)
\end{equation*}
\item{$\bullet$ }  if $f=1_{(a,b]}$ then $(\Gamma ^{(1)}_{H} 1_{(a,b]} )(s)=K_{H}(b,s)-K_{H}(a,s)$.
\item{$\bullet$ } If $f_{i} \in \left| {\cal{H}}\right|$ ($i=1,..,n$) then
\begin{equation}\label{pg3}
\Gamma ^{(n)}_{H}(f_{1}\otimes \ldots \otimes f_{n}) = \Gamma
^{(1)}_{H}f_{1} \otimes \ldots \otimes\Gamma ^{(1)}_{H} f_{n}.
\end{equation}
\end{description}

Let $f\in \left| {\cal{H}}\right| ^{\otimes n}$. Then we define the multiple Wiener-It\^o integral of $f$ with respect to the fractional Brownian motion by
\begin{equation}
\label{def-int}
I_{n}^{H}(f)= I_{n}( \Gamma ^{(n)}f)
\end{equation}
where $I_{n}$ denotes the standard Wiener-It\^o integral with respect to the Wiener process as defined above. Note that $\Gamma ^{(n)}f\in L^{2}([0,1] ^{n})$.

\section{Approximation of multiple fractional Wiener integrals }

Let us introduce some notation. We set
\begin{equation}
\label{eta}
\eta _{\eps } (t) = \int_{0}^{t} K_{H}(t,s) \theta _{\eps } (s)ds, \hskip0.5cm t\in [0,1]
\end{equation}
where $\theta _{\eps }$ is such that $\int_{0}^{t} \theta _{\eps }(s)ds $ converges weakly in the topology of the space ${\cal{C}}_{0}([0,1])$ to the standard Brownian motion.

\begin{lemma}\label{lemaBJT}
Let $\theta_{\eps}$ be either the Kac-Stroock kernels or the Donsker
kernels.  Then the  family of processes $\eta _{\eps}$ converges
weakly in ${\cal{C}}_{0}([0,1])$ as $\eps \to 0$ to the
fractional Brownian motion $B^{H}$ for any $H\in (0,1)$.
\end{lemma}
{\bf Proof: } It has been proved in  \cite{BJT1}, Proposition 2.1. \qed

\vskip0.5cm

Denote, for every $\eps >0$
\begin{equation}\label{diag}
g_{\eps}(x_{1}, \ldots , x_{n})= \prod_{i,j=1; i\not= j}^{n}
1_{\{\vert x_{i}- x_{j} \vert >\eps \}}
\end{equation}
and
\begin{equation}
\label{ieta} I_{n _{\eps}}(f)_{t}= \int_{0}^{1}\ldots \int_{0}^{1}
\left( \Gamma ^{(n)}_{H} f1_{[0,t]}^{\otimes n}\right)(x_{1}, \dots,
x_{n})\left( \prod_{i=1}^{n}\theta _{\eps }(x_{i}) \right)g_{\eps}
(x_{1},\ldots , x_{n}) dx_{1}\ldots dx_{n}.
\end{equation}

\begin{remark}
Note that it follows from a result in \cite{Sa-To} that, if $f\in
L^{q}([0,1])$ for some $q>\frac1H$ then the function $t\to I_{n
_{\eps}}(f)_{t}$ is continuous. Indeed, for every $s<t$
\begin{eqnarray*}
&&\left| I_{n _{\eps}}(f)_{t}-I_{n _{\eps}}(f)_{s}\right| \\
&\leq & \sup _{0\leq r \leq 1}\vert \theta _{\eps} (r) \vert^n \int_{[0,1]^{\otimes n}} \left| (\Gamma _{H} ^{(n)} f1_{[0,t] }^{\otimes n} ) (x_{1},\dots, x_{n}) -(\Gamma _{H} ^{(n)} f1_{[0,s] }^{\otimes n} ) (x_{1},\dots, x_{n})\right| dx_{1}\cdots dx_{n}\\
&\leq & \sup _{0\leq r \leq 1}\vert \theta _{\eps} (r) \vert^n \left( \int_{[0,1]^{\otimes n}} \left| (\Gamma _{H} ^{(n)} f1_{[0,t] }^{\otimes n} ) (x_{1},\dots, x_{n}) -(\Gamma _{H} ^{(n)} f1_{[0,s] }^{\otimes n} ) (x_{1},\dots, x_{n})\right|^{2} dx_{1}\cdots dx_{n}\right) ^{\frac{1}{2}}\\
&=& \sup _{0\leq r \leq 1}\vert \theta _{\eps} (r)  \vert^n  \left( E\left| I_{n}^{H} (f1_{[0,t]}^{\otimes n})-I_{n}^{H} (f1_{[0,s]}^{\otimes n})\right| ^{2}\right) ^{\frac{1}{2}}\\
&\leq & C_{H,n} \sup _{0\leq r \leq 1}\vert \theta _{\eps} (r)
\vert^n\vert t-s\vert ^{H-\frac{1}{q} },
\end{eqnarray*}
where for the last inequality we used Theorem 3.2 in \cite{Sa-To}.
\end{remark}

We first prove the following result.

\begin{lemma}\label{app1}
Let $f$ be a simple function of the form
\begin{equation}\label{step}
f(x_{1}, \dots, x_{n})= \sum_{k=1}^{m} \alpha _{k} 1_{\Delta
_{k}}(x_{1}, \dots, x_{n} )
\end{equation}
where $m\in \mathbb{N}$, $\alpha _{k}\in \mathbb{R}$ for every
$k=1,\dots,$m  and $\Delta _{k}= (a_{k}^{1}, b_{k}^{1}]\times
\cdots\times (a_{k}^{n}, b_{k}^{n}]$ such that  for every
$k=1,\dots,m$, $(a_{k}^{i}, b_{k}^{i}]$ are disjoint intervals
($i=1,\dots,n$). Then the finite dimensional distributions of the
process $(Y^{\eps }(f)_{t})_{t\in [0,1]}$ given by
\begin{equation*}
Y^{\eps }(f)_{t}:=\int_{0}^{1}\ldots \int_{0}^{1} \left( \Gamma
^{(n)}_{H} f1_{[0,t]}^{\otimes n}\right)(x_{1}, \dots, x_{n})\left(
\prod_{i=1}^{n}\theta _{\eps }(x_{i}) \right) dx_{1}\ldots dx_{n}\,,
\end{equation*}
converge as $\eps \to 0$ to the finite dimensional distributions of
\begin{eqnarray*}
&&\left( \sum_{k=1}^{m} \alpha _{k}  I_{1}^{H}(1_{(a_{k}^{1}, b_{k}^{1}]}1_{[0,t]}) \cdots I_{1}^{H}(1_{(a_{k}^{n}, b_{k}^{n}]}1_{[0,t]}) \right)  _{t\in [0,1]}\\
&=& \left( \sum_{k=1}^{m} \alpha _{k}  \left( B^{H}_{b_{k}^{1}\wedge
t} -B^{H}_{a_{k}^{1}\wedge t}\right) \cdots \left(
B^{H}_{b_{k}^{n}\wedge t} -B^{H}_{a_{k}^{n}\wedge t}\right)\right)
_{t\in [0,1]}.
\end{eqnarray*}
\end{lemma}
{\bf Proof: } We have, by using the property (\ref{pg3}) of the operator $\Gamma _{H}^{(n)}$,
\begin{eqnarray*}
Y^{\eps }(f)_{t}&=&
 \int_{0}^{1}\ldots \int_{0}^{1} \left( \Gamma ^{(n)}_{H} f1_{[0,t]}^{\otimes n}\right)(x_{1}, \dots, x_{n})\left( \prod_{i=1}^{n}\theta _{\eps }(x_{i}) \right) dx_{1}\ldots dx_{n}\\
 &=&\sum_{k=1}^{m}\alpha _{k} \int_{0}^{1}\ldots \int_{0}^{1}\left( \Gamma ^{(n)}_{H} 1_{ (a_{k}^{1}, b_{k}^{1}]\times \cdots\times (a_{k}^{n}, b_{k}^{n}] }1_{[0,t]}^{\otimes n}\right) (x_{1}, \ldots , x_{n})\left(\prod_{i=1}^{n}\theta _{\eps }(x_{i})\right)  dx_{1}\ldots dx_{n}\\
 &=&\sum_{k=1}^{m}\alpha _{k}\prod_{i=1}^{n} \int_{0}^{1}\left( \Gamma ^{(1)}_{H} 1_{ (a_{k}^{i}, b_{k}^{i}]}1_{[0,t]}\right)(x_{i})\theta_{\eps}(x_{i})dx_{i}\\
 &=&\sum_{k=1}^{m}\alpha _{k}\prod_{i=1}^{n} \int_{0}^{1}\left(  K_{H}(b_{k}^{i}\wedge t, x_{i})-K_{H}(a_{k}^{i}\wedge t, x_{i})\right)\theta_{\eps}(x_{i})dx_{i}\\
 &=&
 \sum_{k=1}^{m}\alpha _{k}\prod_{i=1}^{n}(\eta_{\eps} (b_{k}^{i}\wedge t) -\eta_{\eps}(a_{k}^{i}\wedge t)).
 \end{eqnarray*}
 and then for every fixed $t_{1}, \dots, t_{r}\in [0,1]$ the vector $(Y^{{\eps}}(f)_{t_{1}}, \ldots ,  Y^{{\eps}}(f)_{t_{r}} )$ converges as in the statement  because by Lemma 1  $\eta _{\eps} $ converges weakly to the fractional Brownian motion. \qed

 \vskip0.2cm

 \begin{remark}
 Let $f$ be a simple function. It can be seen that (here $\partial_{1}K _{H}$ denotes the partial derivative of $K_{H}$ with respect to the first variable)
 \begin{eqnarray*}
&& Y^{\eps }(f)_{t}\\&=&\int_{[0,t] ^{n}}dx_{1}\cdots dx_{n}\left( \int_{x_{1}}^{t}\ldots  \int_{x_{n}}^{t}\partial_{1}K_{H}(t_{1}, x_{1})\ldots \partial_{1}K_{H}(t_{n}, x_{n})f(t_{1},\dots,t_{n}) dt_{1}\cdots dt_{n}\right)\left(\prod_{i=1}^n\theta_{\eps}(x_i)\right) \\
 &=& \int_{[0,t] ^{n}}f(t_{1},\dots,t_{n})d\eta_{\eps} (t_{1}) \cdots d\eta_{\eps}(t_{n}) .
 \end{eqnarray*}
 Therefore $Y^{\eps }(f)$ coincides with  $I^{3}_{n_{\eps}}$ defined by (\ref{i3}).
 \end{remark}
 Note also that in the case of multiple Wiener-It\^o integrals ($H=\frac{1}{2}$) the random variable $Y^{\eps}(f)_{t} $ coincides with $I_{n_{\eps}}(f)_{t}$ for $\eps$ small enough if $f$ is a simple function.

 \subsection{The case $n=2$ }

 Let us consider first the case of a multiple integral in the second Wiener chaos. Suppose that $f$ is a simple function of two variables of the form
\begin{equation*}
f(x, y)= \sum_{k=1}^{m} \alpha _{k} 1_{(a_{k} ^{1}, b_{k}^{1}]}(x) 1_{(a_{k} ^{2}, b_{k}^{2}]}(y)
\end{equation*}
where for every $k$,  $(a_{k} ^{1}, b_{k}^{1}]$ and $(a_{k} ^{2}, b_{k}^{2}]$ are disjoint intervals. In this case, by using the product formula  for multiple stochastic integrals (\ref{prod}),  the multiple integral of $f$ with respect to $B^{H}$ can be expressed as
\begin{eqnarray}
I_{2}^{H}(f1_{[0,t]}^{\otimes 2})&=& \sum_{k=1}^{m}\alpha _{k} I_{2}^{H} \left(  1_{(a_{k} ^{1}, b_{k}^{1}]} 1_{(a_{k} ^{2}, b_{k}^{2}]}1_{[0,t]}^{\otimes 2} \right)\nonumber \\
&=&  \sum_{k=1}^{m}\alpha _{k} \left( B^{H} _{b_{k}^{1}\wedge t}- B^{H} _{a_{k}^{1}\wedge t}\right) \left( B^{H} _{b_{k}^{2}\wedge t}- B^{H} _{a_{k}^{2}\wedge t}\right) \nonumber \\
&& -\sum_{k=1}^{m}\alpha _{k} \langle 1_{(a_{k} ^{1}, b_{k}^{1}]}1_{[0,t]}, 1_{(a_{k} ^{2}, b_{k}^{2}]}1_{[0,t]}\rangle _{{\cal{H}}} .\label{d2}
\end{eqnarray}
The main difference with respect  to the case of the standard Brownian motion is given by the fact that the scalar product in ${\cal{H}}$ of two indicator functions of disjoint intervals is not zero anymore.

Let us show that the sequence
\begin{equation}\label{I2e}
I_{2 _{\eps}}(f)_{t}= \int_{0}^{1}\int_{0}^{1} \left( \Gamma
^{(2)}_{H} f1_{[0,t]}^{\otimes 2}\right)(x_{1}, x_{2})\left(
\prod_{i=1}^{2}\theta _{\eps }(x_{i}) \right)1_{\{\vert
x_{1}-x_{2}\vert >\eps \}}  dx_{1} dx_{2}
\end{equation}
converges in the sense of finite dimensional distributions to the
process $I_2^{H}(f1_{[0,t]}^{\otimes 2} )$. We can write
\begin{eqnarray}
I_{2 _{\eps}}(f)_{t}&=& \int_{0}^{1}\int_{0}^{1} \left( \Gamma ^{(2)}_{H} f1_{[0,t]}^{\otimes 2}\right)(x_{1}, x_{2})\left( \prod_{i=1}^{2}\theta _{\eps }(x_{i}) \right)  dx_{1} dx_{2}\nonumber\\
&&- \int_{0}^{1}\int_{0}^{1} \left( \Gamma ^{(2)}_{H}
f1_{[0,t]}^{\otimes 2}\right)(x_{1}, x_{2})\left(
\prod_{i=1}^{2}\theta _{\eps }(x_{i}) \right)1_{\{\vert
x_{1}-x_{2}\vert <\eps \}}  dx_{1} dx_{2}.\label{d1}
\end{eqnarray}

Note that, using the properties of the transfer operator $\Gamma _{H} ^{(2)}$, the first term can be written as
\begin{eqnarray*}
&&\sum_{k}\alpha _{k} \int_{0} ^{1} \Gamma _{H} ^{(1) } \left( 1_{(a_{k}^{1}, b_{k}^{1}]\cap [0,t] }\right) (x_{1})
\theta _{\eps }(x_{1} ) dx_{1}   \int_{0} ^{1} \Gamma _{H} ^{(1) } \left( 1_{(a_{k}^{2}, b_{k}^{2}]\cap [0,t] }\right) (x_{2}) \theta _{\eps }(x_{2} ) dx_{2}\\
&=& \sum_{k}\alpha _{k} (\eta _{\eps } (b_{k}^{1} \wedge t)  -\eta
_{\eps } (a_{k}^{1} \wedge t)) (\eta _{\eps } (b_{k}^{2} \wedge t)
-\eta _{\eps } (a_{k}^{2} \wedge t))
\end{eqnarray*}
and by Lemma 2,  its finite dimensional distributions  converge to  those of the stochastic process
$$\sum_{k=1}^{m}\alpha _{k} \left( B^{H} _{b_{k}^{1}\wedge t}- B^{H} _{a_{k}^{1}\wedge t}\right) \left( B^{H} _{b_{k}^{2}\wedge t}- B^{H} _{a_{k}^{2}\wedge t }\right).$$

Next we will  discuss the behavior as $\eps \to 0$ of the second
term. We need the following lemma, which  will play an important
role in the sequel.

\begin{lemma}\label{imp}
Consider two functions $f,g\in L^{2}([0,1] )$ and denote by
\begin{equation*}
Y_{\eps} = \int_{0}^{1}\int_{0}^{1}dx_{1}dx_{2} f(x_{1}) g(x_{2})
\theta _{\eps} (x_{1}) \theta _{\eps }(x_{2}) 1_{\{\vert
x_{1}-x_{2}\vert <\eps\}}
\end{equation*}
where $\theta _{\eps}$ are the Kac-Stroock kernels or the Donsker
kernels. Then
\begin{equation*}
 Y_{\eps} \mathop{\longrightarrow}_{{\eps \to 0}} Y= \int_{0}^{1} f(x)g(x)dx \hskip0.5cm \mbox{ in }
L^{2}(\Omega).
\end{equation*}
\end{lemma}
{\bf Proof: } {\it  The case when  $\theta _{\eps}$ are the
Kac-Stroock kernels }. In this case, $\theta
_{\eps}(x)=\frac{1}{\eps}(-1)^{N(\frac{x}{\eps^2})},$ where
$\{N(t);t\geq0\}$ is a standard Poisson process. We have
\[E\left(Y_{\eps}-Y\right)^2=E\left(Y_{\eps}\right)^2-2YE(Y_{\eps})+E\left(Y\right)^2.\]
We first calculate,
\begin{eqnarray*}E\left(Y_{\eps}\right)&=&\int_0^1\int_0^1dx_{1}dx_{2}f(x_{1}) g(x_{2})\frac{1}{\eps^2}e^{\frac{-2}{\eps^2}\vert x_1-x_2\vert}1_{\{\vert x_{1}-x_{2}\vert
<\eps\}}\\&=&\int_0^1dx_{1}f(x_{1})\int_0^{x_1}dx_{2}
g(x_{2})\frac{1}{\eps^2}e^{\frac{-2}{\eps^2}(x_1-x_2)}1_{\{0<
x_{1}-x_{2}
<\eps\}}\\&&+\int_0^1dx_{2}g(x_{2})\int_0^{x_2}dx_{1}f(x_{1})
\frac{1}{\eps^2}e^{\frac{-2}{\eps^2}(x_2-x_1)}1_{\{ 0<x_{2}-x_{1}
<\eps\}}.
\end{eqnarray*}
Note that
\begin{equation*}
\int_0^{x_1}dx_{2}
g(x_{2})\frac{2}{\eps^2}e^{\frac{-2}{\eps^2}(x_1-x_2)}1_{\{0<
x_{1}-x_{2} <\eps\}} = g\ast \varphi _{\eps }(x_{1} )
\end{equation*}
where $\varphi_{\eps} (z)=1_{(0,\eps )}(z) \frac{2}{\eps ^{2}}
e^{-\frac{2z}{\eps ^{2}}}$ is an approximation of the identity.
Therefore the convolution $g\ast \varphi_{\eps }$ converges to $g$
in $L^{2}([0,1])$ because $g\in L^{2}([0,1])$. We obtain
\[E\left(Y_{\eps}\right)\underset{\eps\rightarrow0}{\longrightarrow} \int_{0}^{1} f(x)g(x)dx=Y\,.\]
On the other hand,
\begin{eqnarray*}&&E\left((Y_{\eps})^2\right)\\&=&\int_{[0,1]^4}dx_{1}dx_{2}dx_{3}dx_{4}f(x_{1})g(x_{2})f(x_{3})g(x_{4})E\left(
\theta _{\eps}(x_1)\ldots\theta _{\eps}(x_4)\right)1_{\{\vert
x_{1}-x_{2}\vert <\eps\}}1_{\{\vert x_{3}-x_{4}\vert <\eps\}}
\\&:=&I_1^{\eps}+I_2^{\eps}
\end{eqnarray*}
with
\begin{eqnarray*}I_1^{\eps}&=&\int_{[0,1]^4}dx_{1}dx_{2}dx_{3}dx_{4}f(x_{1})g(x_{2})f(x_{3})g(x_{4})\frac{1}{\eps^4}e^{\frac{-2}{\eps^2}(
\vert x_2-x_1\vert )}e^{\frac{-2}{\eps^2}( \vert x_4-x_3\vert )}
\\&&\qquad \times1_{\{ \vert x_{2}-x_{1}
\vert <\eps\}}1_{\{ \vert x_{4}-x_{3}\vert <\eps\}}\left( 1_{\{
x_{1}\vee x_{2}<x_{3}\wedge x_{4}\}}+1_{\{ x_{3}\vee
x_{4}<x_{1}\wedge x_{2} \}}\right)
\end{eqnarray*}
and
\begin{eqnarray*}I_2^{\eps}&=&\int_{[0,1]^4}dx_{1}dx_{2}dx_{3}dx_{4}f(x_{1})g(x_{2})f(x_{3})g(x_{4})\frac{1}{\eps^4}e^{\frac{-2}{\eps^2}(
 x_{(2)}-x_{(1)})}e^{\frac{-2}{\eps^2}( x_{(4)}-x_{(3)})}
\\&&\qquad \times1_{\{ \vert x_{2}-x_{1}
\vert <\eps\}}1_{\{ \vert x_{4}-x_{3}\vert <\eps\}} 1_{
A}(x_{1},x_{2}, x_{3}, x_{4})
\end{eqnarray*}
where we denoted by $A=\{\{ x_{1}\vee x_{2}<x_{3}\wedge
x_{4}\}\cup\{x_{3}\vee x_{4}<x_{1}\wedge x_{2}\}\}^C$. We begin
studying the convergence of the term $I_2^{\eps}$. In the set $A$,
there are 16 possible orders for the variables $x_1,x_2,x_3,x_4$. We
will make the calculation for the case $x_1<x_3<x_2<x_4$ but for the
other 15 possible orders we can proceed in a similar way. In this
cas we have,
\begin{eqnarray*}&&\int_{[0,1]^4}dx_{1}dx_{2}dx_{3}dx_{4}f(x_{1})g(x_{2})f(x_{3})g(x_{4})\frac{1}{\eps^4}e^{\frac{-2}{\eps^2}(
x_3-x_1)}e^{\frac{-2}{\eps^2}(
x_4-x_2)}\\&&\qquad\qquad\qquad\qquad\qquad\qquad\times1_{\{
x_{1}<x_{3}<x_{2}<x_{4} \}}1_{\{\vert x_{2}-x_{1}\vert
<\eps\}}1_{\{\vert
x_{4}-x_{3}\vert <\eps\}}\\
&\leq& \frac12\left(
\int_{[0,1]^2}dx_{1}dx_{2}\left[f(x_{1})g(x_{2})\right]^21_{\{\vert
x_{2}-x_{1}\vert
<\eps\}}\right. \\
&& \times \left(
\int_{[0,1]^2}dx_{3}dx_{4}\frac{1}{\eps^4}e^{\frac{-2}{\eps^2}(
x_3-x_1)}e^{\frac{-2}{\eps^2}( x_4-x_2)}1_{\{
x_{1}<x_{3}\}}1_{\{x_{2}<x_{4} \}}1_{\{\vert x_{4}-x_{3}\vert
<\eps\}}\right) \\
&&+
\int_{[0,1]^2}dx_{3}dx_{4}\left[f(x_{3})g(x_{4})\right]^21_{\{\vert
x_{4}-x_{3}\vert <\eps\}} \\
&&\left. \times \left( \int_{[0,1]^2}dx_{1}dx_{2}
\frac{1}{\eps^4}e^{\frac{-2}{\eps^2}( x_3-x_1)}e^{\frac{-2}{\eps^2}(
x_4-x_2)}1_{\{ x_{1}<x_{3}\}}1_{\{x_{2}<x_{4}\}}1_{\{\vert
x_{2}-x_{1}\vert <\eps\}} \right)\right).
\end{eqnarray*}
When we integrate the first integral with respect to $dx_{3}dx_{4}$
and the second integral with respect to $dx_{1}dx_{2}$ we obtain
that the last expression can be bounded by
\begin{eqnarray*}
&&
C\int_{[0,1]^2}dx_{1}dx_{2}\left[f(x_{1})g(x_{2})\right]^21_{\{\vert
x_{2}-x_{1}\vert <\eps\}}.
\end{eqnarray*}

Proceeding in a similar way for the other 15 possible orders we
obtain that
\begin{eqnarray*}
I_{2}^{\eps } &\leq&
C\int_{[0,1]^2}dx_{1}dx_{2}\left[f(x_{1})g(x_{2})\right]^21_{\{\vert
x_{2}-x_{1}\vert <\eps\}}.
\end{eqnarray*}
This implies that $I_2^{\eps}$ converges to $0$, by using the
dominated convergence theorem.

Let us regard the behavior of the term $I_{1} ^{\eps}$. This term will give the convergence of $E(Y^{\eps})^{2}$.  We have
\begin{eqnarray*}I_{1}^{\eps}&=&8\int_{[0,1]^4}dx_{1}dx_{2}dx_{3}dx_{4}f(x_{1})g(x_{2})f(x_{3})g(x_{4})\frac{1}{\eps^4}e^{\frac{-2}{\eps^2}(
\vert x_2-x_1\vert )}e^{\frac{-2}{\eps^2}( \vert x_4-x_3\vert )}
\\&&\qquad\qquad\qquad\qquad\qquad\qquad\times 1_{\{ \vert x_{2}-x_{1} \vert <\eps\}}1_{\{ \vert x_{4}-x_{3}\vert
<\eps\}} 1_{\{
x_{1}< x_{2}<x_{3}< x_{4}\}}\\
&=&
2\int_0^1dx_{2}g(x_{2})1_{\{x_{2}<x_{4}\}}\int_0^{x_{2}}dx_{1}f(x_{1})\frac{2}{\eps^2}e^{\frac{-2}{\eps^2}(
x_2-x_1)}1_{\{0< x_{2}-x_{1} <\eps\}}\\&&\times
\int_0^1dx_{4}g(x_{4})\int_{x_{2}}^{x_4}dx_{3}f(x_{3})\frac{2}{\eps^2}e^{\frac{-2}{\eps^2}(
x_4-x_3)}1_{\{0< x_{4}-x_{3} <\eps\}}.
\end{eqnarray*}
We obtain that $I_{1}^{\eps}$ converges to
$2\left[\int_{0}^{1}\int_{0}^{1}
f(x)g(x)f(y)g(y)1_{\{x<y\}}dxdy\right]$.  Thus $I_1^{\eps}$
converges to $\left[\int_{0}^{1} f(x)g(x)dx\right]^2=Y^2.$

\vskip0.2cm

{\it  The case when  $\theta _{\eps}$ are the Donsker kernels }. In
this case,
$\theta_{\varepsilon}(x)=\frac{1}{\varepsilon}\sum_{k=1}^{\infty}\xi_k1_{[k-1,k)}(\frac{x}{\varepsilon^2})$
where $(\xi_k)$ is a sequence of independent, identically
distributed random variables satisfying $E(\xi_1)=0$ and
$E(\xi_1^2)=1$ with $E(\xi_1^{2n})<+\infty$. In this case we have
\begin{eqnarray*}E\left(Y_{\eps}\right)&=&\int_0^1\int_0^1dx_{1}dx_{2}f(x_{1}) g(x_{2})E\left(\theta_{\varepsilon}(x_1)\theta_{\varepsilon}(x_2)\right)1_{\{\vert x_{1}-x_{2}\vert
<\eps\}}\\&=&\frac{1}{\varepsilon^2}\sum_{k=1}^{\infty}\int_0^1\int_0^{1}dx_{1}dx_{2}f(x_{1})
g(x_{2})1_{[k-1,k)^2}\left(\frac{x_1}{\varepsilon^2},\frac{x_2}{\varepsilon^2}\right)\\&=&
\int_0^1dx_{1}f(x_{1})\sum_{k=1}^{\infty}1_{[(k-1)\varepsilon^2,k\varepsilon^2)}
\left(x_1\right)\frac{1}{\varepsilon^2}\int_{(k-1)\varepsilon^2}^{k\varepsilon^2}dx_{2}g(x_{2})
\\&=&
\int_0^1dx_{1}f(x_{1})\sum_{k=1}^{\left[\frac{1}{\varepsilon^2}\right]+1}1_{[(k-1)\varepsilon^2,k\varepsilon^2)}
\left(x_1\right)\frac{1}{\varepsilon^2}\int_{(k-1)\varepsilon^2}^{k\varepsilon^2}dx_{2}g(x_{2})\\&=&\int_0^1dx_{1}f(x_{1})G_{\varepsilon}(x_1),
\end{eqnarray*}
where
$$G_{\varepsilon}(x):=\sum_{k=1}^{\left[\frac{1}{\varepsilon^2}\right]+1}1_{[(k-1)\varepsilon^2,k\varepsilon^2)}
\left(x\right)\frac{1}{\varepsilon^2}\int_{(k-1)\varepsilon^2}^{k\varepsilon^2}dyg(y).$$
Fix $x_1\in(0,1)$. Then for every $\varepsilon>0$ close to zero,
there exists an
$k(x_1,\varepsilon)\in\{1,\ldots,\left[\frac{1}{\varepsilon^2}\right]+1\}$
such that $(k(x_1,\varepsilon)-1)\varepsilon^2\leq
x_1<k(x_1,\varepsilon)\varepsilon^2$. Then $0\leq
x_1-(k(x_1,\varepsilon)-1)\varepsilon^2<\varepsilon^2$ and this
implies that $(k(x_1,\varepsilon)-1)\varepsilon^2\rightarrow x_1$ as
$\varepsilon \rightarrow0$. Thus
\[G_{\varepsilon}(x_1)=\frac{1}{\varepsilon^2}\int_{(k(x_1,\varepsilon)-1)\varepsilon^2}^{k(x_1,\varepsilon)\varepsilon^2}dx_{2}g(x_{2})\]
converges to $g(x_1)$ as $\varepsilon \rightarrow0$. Consequently
\begin{eqnarray*}E\left(Y_{\eps}\right)\underset{\varepsilon
\rightarrow0}{\longrightarrow}\int_0^1f(x_{1})g(x_1)dx_{1}.
\end{eqnarray*}
Now, we calculate $E\left(Y_{\eps}^2\right)$.
We have that,
\begin{eqnarray*}&&E\left(
\theta _{\eps}(x_1)\ldots\theta
_{\eps}(x_4)\right)\\
&=&\frac{1}{\varepsilon^4}\sum_{k\neq
j=1}^{\infty}1_{[(k-1)\varepsilon^2,k\varepsilon^2)}(x_1,x_2)1_{[(j-1)\varepsilon^2,j\varepsilon^2)}(x_3,x_4)
\\&&+\frac{1}{\varepsilon^4}\sum_{k\neq
j=1}^{\infty}1_{[(k-1)\varepsilon^2,k\varepsilon^2)^2}(x_1,x_3)
1_{[(j-1)\varepsilon^2,j\varepsilon^2)^2}(x_2,x_4)
\\&&+\frac{1}{\varepsilon^4}\sum_{k\neq
j=1}^{\infty}1_{[(k-1)\varepsilon^2,k\varepsilon^2)^2}(x_1,x_4)
1_{[(j-1)\varepsilon^2,j\varepsilon^2)^2}(x_2,x_3)
\\&&+\frac{E(\xi_1^4)}{\varepsilon^4}\sum_{k=1}^{\infty}
1_{[(k-1)\varepsilon^2,k\varepsilon^2)^4}(x_1,x_2,x_3,x_4)\\&:=&G_{\varepsilon}^1(x_1,x_2,x_3,x_4)+G_{\varepsilon}^2(x_1,x_2,x_3,x_4)+G_{\varepsilon}^3(x_1,x_2,x_3,x_4)
+G_{\varepsilon}^4(x_1,x_2,x_3,x_4).
\end{eqnarray*}
Thus
\begin{eqnarray*}E\left((Y_{\eps})^2\right)&=&\int_{[0,1]^4}dx_{1}dx_{2}dx_{3}dx_{4}f(x_{1})g(x_{2})f(x_{3})g(x_{4})G_{\varepsilon}^1(x_1,x_2,x_3,x_4)1_{\{\vert
x_{1}-x_{2}\vert <\eps\}}1_{\{\vert x_{3}-x_{4}\vert
<\eps\}}\\&&+\int_{[0,1]^4}dx_{1}dx_{2}dx_{3}dx_{4}f(x_{1})g(x_{2})f(x_{3})g(x_{4})G_{\varepsilon}^2(x_1,x_2,x_3,x_4)1_{\{\vert
x_{1}-x_{2}\vert <\eps\}}1_{\{\vert x_{3}-x_{4}\vert
<\eps\}}\\&&+\int_{[0,1]^4}dx_{1}dx_{2}dx_{3}dx_{4}f(x_{1})g(x_{2})f(x_{3})g(x_{4})G_{\varepsilon}^3(x_1,x_2,x_3,x_4)1_{\{\vert
x_{1}-x_{2}\vert <\eps\}}1_{\{\vert x_{3}-x_{4}\vert
<\eps\}}\\&&+\int_{[0,1]^4}dx_{1}dx_{2}dx_{3}dx_{4}f(x_{1})g(x_{2})f(x_{3})g(x_{4})G_{\varepsilon}^4(x_1,x_2,x_3,x_4)1_{\{\vert
x_{1}-x_{2}\vert <\eps\}}1_{\{\vert x_{3}-x_{4}\vert <\eps\}}
\\&:=&J_{\eps}^1+J_{\eps}^2+J_{\eps}^3+J_{\eps}^4.
\end{eqnarray*}
{\it The convergence of  $J_{\eps}^1$:} Fix $x_1$ and $x_3$ in
$[0,1]$. Then for every $\varepsilon>0$ close to zero, there exist
$k(x_1,\varepsilon)\in\{1,\ldots,\left[\frac{1}{\varepsilon^2}\right]+1\}$
and
$j(x_3,\varepsilon)\in\{1,\ldots,\left[\frac{1}{\varepsilon^2}\right]+1\}$
such that $k(x_1,\varepsilon)\neq j(x_3,\varepsilon) $,
$(k(x_1,\varepsilon)-1)\varepsilon^2\leq
x_1<k(x_1,\varepsilon)\varepsilon^2$ and
$(j(x_3,\varepsilon)-1)\varepsilon^2\leq
x_3<j(x_3,\varepsilon)\varepsilon^2$,  this implies that
$(k(x_1,\varepsilon)-1)\varepsilon^2\rightarrow x_1$ and
$(j(x_3,\varepsilon)-1)\varepsilon^2\rightarrow x_3$ as $\varepsilon
\rightarrow0$. Then we can write
\[\int_{[0,1]^2}dx_{2}dx_{4}g(x_{2})g(x_{4})G_{\varepsilon}^1(x_1,x_2,x_3,x_4)=
\frac{1}{\varepsilon^2}\int_{(k(x_1,\varepsilon)-1)\varepsilon^2}^{k(x_1,\varepsilon)\varepsilon^2}dx_{2}g(x_{2})\times
\frac{1}{\varepsilon^2}\int_{(j(x_3,\varepsilon)-1)\varepsilon^2}^{j(x_3,\varepsilon)\varepsilon^2}dx_{4}g(x_{4}).\]
Moreover, this term converges to $g(x_{1})g(x_{3})$. We conclude
that
\begin{eqnarray}J_{\eps}^1\underset{\varepsilon
\rightarrow0}{\longrightarrow}\int_0^1\int_0^1dx_{1}dx_{3}f(x_{1})g(x_1)f(x_{3})g(x_3)=Y^2.
\end{eqnarray}
\vskip0.2cm {\it The convergence of  $J_{\eps}^2$ and $J_{\eps}^3$
:} Fix $x_1$ and $x_2$ in $[0,1]$. Then for every $\varepsilon>0$
close to zero, there exist
$k(x_1,\varepsilon)\in\{1,\ldots,\left[\frac{1}{\varepsilon^2}\right]+1\}$
and
$j(x_2,\varepsilon)\in\{1,\ldots,\left[\frac{1}{\varepsilon^2}\right]+1\}$
such that $k(x_1,\varepsilon)\neq j(x_2,\varepsilon) $,
$(k(x_1,\varepsilon)-1)\varepsilon^2\leq
x_1<k(x_1,\varepsilon)\varepsilon^2$ and
$(j(x_2,\varepsilon)-1)\varepsilon^2\leq
x_2<j(x_2,\varepsilon)\varepsilon^2$,  this implies that
$(k(x_1,\varepsilon)-1)\varepsilon^2\rightarrow x_1$ and
$(j(x_2,\varepsilon)-1)\varepsilon^2\rightarrow x_2$ as $\varepsilon
\rightarrow0$. Hence
\begin{eqnarray*}&&\int_{[0,1]^2}dx_{3}dx_{4}g(x_{3})g(x_{4})G_{\varepsilon}^2(x_1,x_2,x_3,x_4)1_{\{\vert
x_{1}-x_{2}\vert <\eps\}}1_{\{\vert x_{3}-x_{4}\vert
<\eps\}}\\&&\leq 1_{\{\vert x_{1}-x_{2}\vert <\eps\}}
\frac{1}{\varepsilon^2}\int_{(k(x_1,\varepsilon)-1)\varepsilon^2}^{k(x_1,\varepsilon)\varepsilon^2}dx_{3}g(x_{3})\times
\frac{1}{\varepsilon^2}\int_{(j(x_2,\varepsilon)-1)\varepsilon^2}^{j(x_2,\varepsilon)\varepsilon^2}dx_{4}g(x_{4}).
\end{eqnarray*} This last term converges to $g(x_{1})g(x_{2})1_{\{
x_{1}=x_{2}\}}$, this implies that
\begin{eqnarray}J_{\eps}^2\underset{\varepsilon
\rightarrow0}{\longrightarrow}0.
\end{eqnarray}
In the same way, we obtain that
\begin{eqnarray}J_{\eps}^3\underset{\varepsilon
\rightarrow0}{\longrightarrow}0.
\end{eqnarray}
{\it The convergence of  $J_{\eps}^4$: }Fix $x_1$  in $[0,1]$. Then for
every $\varepsilon>0$ close to zero, there exist
$k(x_1,\varepsilon)\in\{1,\ldots,\left[\frac{1}{\varepsilon^2}\right]+1\}$
such that $(k(x_1,\varepsilon)-1)\varepsilon^2\leq
x_1<k(x_1,\varepsilon)\varepsilon^2$,  this implies that
$(k(x_1,\varepsilon)-1)\varepsilon^2\rightarrow x_1$ as $\varepsilon
\rightarrow0$. Then we can write
\begin{eqnarray*}&&\int_{[0,1]^3}dx_{2}dx_{3}dx_{4}g(x_{2})f(x_{3})g(x_{4})G_{\varepsilon}^4(x_1,x_2,x_3,x_4)=
\frac{E(\xi_1^4)}{\varepsilon^2}\int_{(k(x_1,\varepsilon)-1)\varepsilon^2}^{k(x_1,\varepsilon)\varepsilon^2}dx_{2}g(x_{2})\\&&\qquad\qquad\qquad\qquad\qquad\times
\frac{1}{\varepsilon^2}\int_{(k(x_1,\varepsilon)-1)\varepsilon^2}^{k(x_1,\varepsilon)\varepsilon^2}dx_{3}f(x_{3})\times
\int_{(k(x_1,\varepsilon)-1)\varepsilon^2}^{k(x_1,\varepsilon)\varepsilon^2}dx_{4}g(x_{4}).\end{eqnarray*}
The last term converges to zero, thus
\begin{eqnarray}J_{\eps}^4\underset{\varepsilon
\rightarrow0}{\longrightarrow}0.
\end{eqnarray}
Consequently, by combining the above convergences we obtain that
\begin{eqnarray}E\left((Y_{\eps})^2\right)\underset{\varepsilon
\rightarrow0}{\longrightarrow}Y^2.
\end{eqnarray}
 \qed

 \vskip0.3cm

We will also need the following lemma.

\begin{lemma}
\label{lema4}
Let us consider a family of stochastic processes $(X^{\eps}) _{t\in [0,1]}$ converging as $\eps \to 0$ to $(X_{t})_{t\in [0,1]}$ in the sense of finite dimensional distributions and a family of stochastic processes $(Y^{\eps}) _{t\in [0,1]}$ such that  for every $t\in [0,1]$ the sequence of random variables $Y^{\eps}_{t}$ converges, as $\eps \to 0$ to $Y_{t}$  in $L^{2}(\Omega)$ where $Y_{t}$ is a constant for every $t$. Then $X^{\eps} + Y^{\eps}$ converges to $X+Y$ in the sense of finite dimensional distributions.

\end{lemma}
{\bf Proof: } Fix $t_{1}, \dots, t_{r}\in [0,1]$ and let us show
that the vector
\begin{equation*}
\left( X^{\eps } _{t_{1}}+ Y^{\eps } _{t_{1}}, \dots,  X^{\eps }
_{t_{r}}+ Y^{\eps } _{t_{r}}  \right)
\end{equation*}
converges in law to the vector
\begin{equation*}\left( X _{t_{1}}+ Y _{t_{1}}, \dots,  X_{t_{r}}+ Y _{t_{r}}  \right) .
\end{equation*} Take $g\in\mathcal C_b^1(\mathbb R^r)$, then
\begin{eqnarray*}
&&|E\left( g(X^{\eps } _{t_{1}}+ Y^{\eps } _{t_{1}}, \dots,  X^{\eps } _{t_{r}}+ Y^{\eps } _{t_{r}})\right) -E\left( g(X _{t_{1}}+ Y _{t_{1}}, \dots,  X _{t_{r}}+ Y _{t_{r}})\right)|\\
&\leq&|E\left( g(X^{\eps } _{t_{1}}+ Y^{\eps } _{t_{1}}, \dots,  X^{\eps } _{t_{r}}+ Y^{\eps } _{t_{r}})\right) -E\left( g(X ^{\eps}_{t_{1}}+ Y _{t_{1}}, \dots,  X^{\eps} _{t_{r}}+ Y _{t_{r}})\right)|\\
&&+ |E\left( g(X ^{\eps}_{t_{1}}+ Y _{t_{1}}, \dots,  X^{\eps}
_{t_{r}}+ Y _{t_{r}})\right)-E\left( g(X _{t_{1}}+ Y _{t_{1}},
\dots, X _{t_{r}}+ Y _{t_{r}})\right)|.
\end{eqnarray*}
The first term converges to zero due to the $L^{2}$ convergence of $Y^{\eps}$ to $Y$ since
\begin{eqnarray*}
&&|E\left( g(X^{\eps } _{t_{1}}+ Y^{\eps } _{t_{1}}, \dots,  X^{\eps } _{t_{r}}+ Y^{\eps } _{t_{r}})\right) -E\left( g(X ^{\eps}_{t_{1}}+ Y _{t_{1}}, \dots,  X^{\eps} _{t_{r}}+ Y _{t_{r}})\right)|\\
&& \leq K E\left[ (Y^{\eps }_{t_{1}} -Y_{t_{1}})^{2}+\cdots+
(Y^{\eps }_{t_{r}} -Y_{t_{r}})^{2}\right]^{\frac12}
\end{eqnarray*}
and the second one converges to zero as $\eps \to 0$ because, by Slutsky's theorem,  $X^{\eps}+Y$ converges to $X+Y$ in the sense of finite dimensional distributions. \qed

\vskip0.3cm

We obtain the following result:

\begin{prop}
Let $f\in {\cal{S}}_2$ and let $I_{2_{\eps}}(f)_{t}$ be given by
(\ref{I2e}). Then $ \left( I_{2_{\eps }} (f)_{t}\right) _{t\in
[0,1]}$ converges as  $\eps \to 0$ in the sense of finite
dimensional distributions to the process $\left( I_{2}^{H}(f1_{[0,t]
}^{\otimes 2} )\right) _{t\in [0,1]}$.
\end{prop}
{\bf Proof: } Recall the expressions (\ref{d2}) and (\ref{d1}) of
$I_{2}^{H}(f1_{[0,t] }^{\otimes 2})$ and $I_{2_{\eps}}(f)_{t}$. By
Lemma 2 the first term in (\ref{d2}) converges in the sense of
finite dimensional distributions to the first term in (\ref{d1}) and
applying Lemma \ref{imp} for $f= \Gamma ^{(1)}_{H} 1_{(a_{k}^{1},
b_{k}^{1}]}1_{[0,t]} $ and $g =\Gamma ^{(1)}_{H} 1_{(a_{k}^{2},
b_{k}^{2}]}1_{[0,t]} $ we obtain that the term
$$
\int_{0}^{1}\int_{0}^{1} \left( \Gamma ^{(2)}_{H}
f1_{[0,t]}^{\otimes 2}\right)(x_{1}, x_{2})\left(
\prod_{i=1}^{2}\theta _{\eps }(x_{i}) \right)1_{\{\vert
x_{1}-x_{2}\vert <\eps\} }  dx_{1} dx_{2}
$$ converges in $L^{2} (\Omega)$  for every $t\in [0,1]$ to
\begin{equation*}
\sum_{k}\alpha _{k}   \int _{0}^{1} \Gamma _{H}^{(1)}
(1_{(a_{k}^{1}, b_{k}^{1}]} 1_{[0,t] } )(x) \Gamma _{H}^{(1)}
(1_{(a_{k}^{2}, b_{k}^{2}]} 1_{[0,t] } ) (x)dx=
 \sum_{k}\alpha _{k} \langle 1_{(a_{k}^{1}, b_{k}^{1}]} 1_{[0,t] }, 1_{(a_{k}^{2}, b_{k}^{2}]} 1_{[0,t] }\rangle _{{\cal{H}}}.
\end{equation*}
The above Lemma 4 gives the conclusion. \qed

\subsection{The case $n\geq 3$}

In the case of multiple integrals of order $n\geq 3$, the structure of $I_{n_{\eps}}(f) _{t}$ is more complex because of the appearance of all diagonals.  The first step is to express the multiple integral of a tensor product of one-variable functions.

\begin{lemma}
Let $f_{1}, \ldots , f_{n} \in {\cal{H}}$. Then
\begin{eqnarray}
\label{tensor}
&& I_{n}^{H}(f_{1}\otimes f_{2} \ldots \otimes f_{n})\nonumber\\
&=&\prod_{i=1}^n I_1^H(f_i)\\
&&+\sum_{l=1} ^{[n/2] } (-1) ^{l}
\sum_{\footnotesize\begin{array}{c}k_{1}, \dots, k_{2l}=1;\\
k_{j}\,\, distinct\end{array}} ^{n} \left( \prod _{u\in
\{1,\dots,n\}\setminus \{k_{1},\dots, k_{2l}\} } I_{1}^{H} (f_{u})
\right)
 \langle f_{k_{1}}, f_{k_{2}}\rangle _{{\cal{H}}}\cdots  \langle f_{k_{2l-1}}, f_{k_{2l}}\rangle _{{\cal{H}}}.\nonumber
\end{eqnarray}
\end{lemma}
{\bf Proof: } We will prove the result by induction. For $n=1,2$ it
is trivial. Let us show how it works for $n=3$ because it is useful
to understand the general case. We have, using (\ref{pg3}),
(\ref{def-int}) and the product formula for multiple integrals
(\ref{prod})
\begin{eqnarray*}
I_{3}^{H}\left( f_{1} \otimes f_{2} \otimes f_{3} \right) &=& I_{3}\left( \Gamma^{(3)}_H(f_{1} \otimes f_{2} \otimes f_{3} )\right)\\
&=&I_{3}\left( \Gamma^{(1)}_H(f_{1}) \otimes \Gamma^{(1)}_H(f_{2}) \otimes \Gamma^{(1)}_H(f_{3})\right)\\
&=&I_{3}\left( \left(\Gamma^{(1)}_H(f_{1}) \otimes \Gamma^{(1)}_H(f_{2}) \otimes \Gamma^{(1)}_H(f_{3})\right)^{\sim}\right)\\
&=& I_{2} \left(\Gamma^{(1)}_H(f_{1}) \tilde{\otimes }
\Gamma^{(1)}_H(f_{2})\right)
I_{1}\left(\Gamma^{(1)}_H(f_{3})\right)\\&& -2 I_{1} \left(\left(
\Gamma^{(1)}_H(f_{1}) \tilde{\otimes } \Gamma^{(1)}_H(f_{2})\right)
\otimes_ {1} \Gamma^{(1)}_H(f_{3}) \right).
\end{eqnarray*}
 Note that
 \begin{equation*}
 \left(\Gamma^{(1)}_H(f_{1})\tilde{\otimes }\Gamma^{(1)}_H(f_{2}) \right) (t_{1}, t_{2}) =\frac{1}{2} \left( \Gamma^{(1)}_H(f_{1})(t_{1}) \Gamma^{(1)}_H(f_{2})
  (t_{2}) + \Gamma^{(1)}_H(f_{1}) (t_{2}) \Gamma^{(1)}_H(f_{2})(t_{1}) \right)
 \end{equation*}
and thus
\begin{equation*}
\left(\Gamma^{(1)}_H(f_{1}) \tilde{\otimes }
\Gamma^{(1)}_H(f_{2})\right) \otimes_ {1}
\Gamma^{(1)}_H(f_{3})=\frac{1}{2} \left( \langle f_{1}, f_{3}\rangle
_{{\cal{H}}}\Gamma^{(1)}_H(f_{2}) + \langle f_{2}, f_{3}\rangle
_{{\cal{H}}} \Gamma^{(1)}_H(f_{1})\right) .
\end{equation*}
We obtain
\begin{eqnarray*}
I_{3}^{H}\left( f_{1} \otimes f_{2} \otimes f_{3} \right)&=& \left( I_{1}^H(f_{1} ) I_{1}^H(f_{2}) -\langle f_{1}, f_{2}\rangle _{{\cal{H}}}\right) I_{1}^{H} (f_{3}) \\
&&+  \left( \langle f_{1}, f_{3}\rangle _{{\cal{H}}}I_{1}^{H} (f_{2}) + \langle f_{2}, f_{3}\rangle _{{\cal{H}}} I_{1}^{H}(f_{1})\right) \\
&=& I_{1}^{H}(f_{1})I_{1}^{H}(f_{2})I_{1}^{H}(f_{3})\\&&- \left(
\langle f_{1}, f_{3}\rangle _{{\cal{H}}}I_{1}^{H} (f_{2}) + \langle
f_{2}, f_{3}\rangle _{{\cal{H}}} I_{1}^{H}(f_{1})+ \langle f_{1},
f_{2}\rangle _{{\cal{H}}}I_{1}^{H}(f_{3})\right).
\end{eqnarray*}
Concerning the  general case, assume that (\ref{tensor}) holds for
$1,2,\dots, n-1$.  Again by the multiplication formula (\ref{prod}),
\begin{eqnarray*}
&&I_{n}^{H}(f_{1}\otimes \cdots\otimes f_{n}) \\
&=&I_{n-1} \left(\Gamma^{(1)}_H(f_{1})\otimes \cdots\otimes
\Gamma^{(1)}_H(f_{n-1})\right)I_{1}\left(\Gamma^{(1)}_H(f_{n})\right)\\&&
-(n-1) I_{n-2} \left( \left(
\Gamma^{(1)}_H(f_{1})\otimes \cdots\otimes \Gamma^{(1)}_H(f_{n-1})\right)^{\sim }\otimes _{1} \Gamma^{(1)}_H(f_{n})\right) \\
&=&I_1^H(f_n)\cdot\prod_{i=1}^{n-1} I_1^H(f_i)\\&&+
 \sum_{l=1} ^{[(n-1)/2] } (-1) ^{l}
\!\!\!\!\!\!\!\!\!\!\sum_{\footnotesize\begin{array}{c}k_{1},\dots, k_{2l}=1;\\
k_{j}\,\, distinct\end{array}}^{n-1}\!\!\!\!\!
 \left( \prod _{u\in \{1,\dots,n-1\}\setminus \{k_{1},\dots, k_{2l}\} } I_{1}^{H} (f_{u}) I_{1}^{H}(f_{n})\right)
 \langle f_{k_{1}}, f_{k_{2}}\rangle _{{\cal{H}}}\dots  \langle f_{k_{2l-1}}, f_{k_{2l}}\rangle _{{\cal{H}}}\\
 &&- \sum_{i=1}^{n-1} I_{n-2}\left( \left(\Gamma^{(1)}_H(f_{1}) \otimes \cdots\hat{i}\cdots \otimes \Gamma^{(1)}_H(f_{n-1})\right)^{\sim}\right)
 \langle f_{i}, f_{n}\rangle _{{\cal{H}}}
\end{eqnarray*}
and this equal to
\begin{eqnarray*}
&&I_1^H(f_n)\cdot\prod_{i=1}^{n-1} I_1^H(f_i)\\&&+ \sum_{l=1}
^{[(n-1)/2] } (-1) ^{l}
\!\!\!\!\!\!\!\!\!\!\sum_{\footnotesize\begin{array}{c}k_{1}, \dots, k_{2l}=1;\\
k_{j}\,\, distinct\end{array}}^{n-1}\!\!\!\!\! \left( \prod _{u\in
\{1, \dots,n-1\}\setminus \{k_{1},\dots, k_{2l}\} }
 I_{1}^{H} (f_{u})I_{1}^{H}(f_{n}) \right)
 \langle f_{k_{1}}, f_{k_{2}}\rangle _{{\cal{H}}}\dots  \langle f_{k_{2l-1}}, f_{k_{2l}}\rangle _{{\cal{H}}}\\
 && -\sum_{i=1}^{n-1}
\big(\prod_{j=1;\,j\neq i}^{n-1} I_1^H(f_j)\\&&+
 \sum_{l=1} ^{[(n-2)/2] } (-1) ^{l}  \sum_{\footnotesize\begin{array}{c}k_{1}, \dots, k_{2l}=1; k_{j}\not= i;\\ k_{j}\,\, distinct\end{array}} ^{n-1} \left( \prod _{u\in \{1,\dots, \hat{i},\dots,n-1\}\setminus \{k_{1},\dots, k_{2l}\} } I_{1}^{H} (f_{u}) \right)\\
 &&\times
\langle f_{k_{1}}, f_{k_{2}}\rangle _{{\cal{H}}}\dots \langle
f_{k_{2l-1}}, f_{k_{2l}}\rangle _{{\cal{H}}}\big)\langle f_{i} ,
f_{n}\rangle _{{\cal{H}}}
\end{eqnarray*}
and it is not difficult to see that the last quantity is equal to the right side of (\ref{tensor}).  \qed

\vskip0.3cm

The next auxiliary two lemmas will be used in the proof of the main result.

\begin{lemma}\label{lema6}
Suppose that $(X^{\eps} _{t})_{t\in [0,1]}$  is a family of
stochastic processes whose finite dimensional distributions
converges to  the finite dimensional distributions  of a stochastic
processes  $(X_{t})_{t\geq 0}$. Suppose also that
$(Y^{\eps}_{s,t})_{s,t\in [0,1]}$  is a two-parameter stochastic
process such that for every $s,t \in [0,1]$ we have  that $Y^ {\eps
} _{s,t}$ converge in  $L^2(\Omega)$ to $a_{s,t}$, when $\eps$ tends
to $0$, where $a_{s,t}$ is a real constant. Then for every
$t_{1},\dots, t_{r_{1}}\in [0,1]$ and $s_{1}, \dots, s_{r_{2}},
u_{1}, \dots, u_{r_{2}} \in [0,1]$ the vector
\begin{equation*}
(X^{\eps} _{t_{1}}, \dots,X^{\eps } _{t_{r_{1}}},  Y^{\eps }
_{s_{1}, u_{1}}, \dots,  Y^{\eps } _{s_{r_{2}}, u_{r_{2}}})
\end{equation*}
converges weakly to the vector
\begin{equation*}
(X _{t_{1}}, \dots,X _{t_{r_{1}}},  a _{s_{1}, u_{1}}, \dots,  a
_{s_{r_{2}}, u_{r_{2}}}).
\end{equation*}
\end{lemma}
{\bf Proof: } Consider $f\in\mathcal C_b^1(\mathbb R^{r_1+r_2})$. We
have
\begin{eqnarray*}&&|E\left(f(X^{\eps} _{t_{1}}, \dots,X^{\eps } _{t_{r_{1}}},
  Y^{\eps } _{s_{1}, u_{1}}, \dots,  Y^{\eps } _{s_{r_{2}}, u_{r_{2}}})\right)
-E\left(f(X _{t_{1}}, \dots,X _{t_{r_{1}}},  a _{s_{1}, u_{1}},
\dots,  a_{s_{r_{2}}, u_{r_{2}}})\right)|\\&\leq&|E\left(f(X^{\eps}
_{t_{1}}, \dots,X^{\eps } _{t_{r_{1}}},  Y^{\eps } _{s_{1}, u_{1}},
\dots,  Y^{\eps } _{s_{r_{2}}, u_{r_{2}}})\right)
-E\left(f(X^{\eps}_{t_{1}}, \dots,X^{\eps}_{t_{r_{1}}},  a_{s_{1},
u_{1}}, \dots, a_{s_{r_{2}},
u_{r_{2}}})\right)|\\&&+|E\left(f(X^{\eps}_{t_{1}},
\dots,X^{\eps}_{t_{r_{1}}}, a_{s_{1}, u_{1}}, \dots, a_{s_{r_{2}},
u_{r_{2}}})\right) -E\left(f(X _{t_{1}}, \dots,X _{t_{r_{1}}},  a
_{s_{1}, u_{1}}, \dots,  a _{s_{r_{2}}, u_{r_{2}}})\right)|
\\&\leq&K E\left[( Y^{\eps } _{s_{1}, u_{1}}-a _{s_{1}, u_{1}})^2+\cdots+(Y^{\eps } _{s_{r_{2}},
u_{r_{2}}}-a_{s_{r_{2}}, u_{r_{2}}})^2 \right]^{\frac12}
\\&&+|E\left(f(X^{\eps}_{t_{1}}, \dots,X^{\eps}_{t_{r_{1}}}, a_{s_{1}, u_{1}}, \dots,
a_{s_{r_{2}}, u_{r_{2}}})\right) -E\left(f(X _{t_{1}}, \dots,X
_{t_{r_{1}}},  a _{s_{1}, u_{1}}, \dots,  a _{s_{r_{2}},
u_{r_{2}}})\right)|.
\end{eqnarray*}
The first term converges to zero because $Y^ {\eps } _{s,t}$
converge in $L^2(\Omega)$ to $a_{s,t}$ and the second one converge
to 0 because the finite dimensional distributions of $(X^{\eps}
_{t})_{t\in [0,1]}$ converge weakly to those of $(X_{t})_{t\geq 0}$. \qed

\begin{lemma}
\label{convR} Let $f_{i}\in L^{2}([0,1] )$ for $i =1, \dots,n$ and
define
\begin{equation*}
F_{\eps}= \int_{[0,1] ^{n}}f_1(x_{1})\cdots f_n(x_{n})1_{\{\vert
x_{1}-x_{2}\vert <\eps\}}1_{\{\vert x_{1}- x_{3}\vert <\eps\}}
\left(\prod_{j=1}^{n}\theta _{\eps} (x_{j} ) \right)h_{\eps}(x_{1},
\dots, x_{n}) dx_{1}\cdots dx_{n}
\end{equation*}
where $\theta _{\eps}$ are the Kac-Stroock's or the Donsker's
kernels and  we assume that for every $\eps >0$, $\vert h_{\eps}
(x_{1},\dots, x_{n}) \vert \leq 1$ for every $x_{1}, \dots, x_{n}\in
[0,1]$. Then
\begin{eqnarray}\label{convergence of F{epsilon}} E(F_{\eps}^2)\underset{\eps\rightarrow0}{\longrightarrow}0.
\end{eqnarray}
\end{lemma}
{\bf Proof: } {\bf 1)} Assume  that $ \theta _{\eps}$ are
Kac-Stroock kernels. For $\eps>0$ and $x\geq0$, we set
$Q_{\eps}(x)=\frac{1}{\eps^2}e^{\frac{-2x}{\eps^2}}$ and
$f_j:=f_{j-n}$ for every $j=n+1,\ldots,2n$. We introduce some
operators on the set of permutations. $\mathfrak{S}_k$ denote  the
set of permutations on ${1, . . . , k}$. When
$\tau\in\mathfrak{S}_{2m}$ and $\sigma\in\mathfrak{S}_{m}$, we note
$\sigma\star \tau$ the element of $\mathfrak{S}_{2m}$ defined by
\[\sigma\star \tau(2j-1)=\tau(2\sigma(j)-1) \mbox{ and }\sigma\star \tau(2j)=\tau(2\sigma(j)).\]
 We have $id\star\tau=\tau$ and $\sigma'\star(\sigma\star\tau)=(\sigma'\sigma)\star\tau$, hence $\star :\mathfrak{S}_{m}\times\mathfrak{S}_{2m} \rightarrow\mathfrak{S}_{2m}$
 defines  a (right) group action of $\mathfrak{S}_{m}$ on
$\mathfrak{S}_{2m}$. For any $\tau\in\mathfrak{S}_{2m}$, the orbit
of $\tau$ has exactly $m!$ elements. Consequently, the set
$\mathcal{O}$ of the orbits under the group action $\star$ has
$\frac{(2m)!}{m!}$  elements and we have, by denoting $\tau_i$ one
particular element of the orbit $o_i=o(\tau_i)\in\mathcal{O}$ : for
$r_1,\ldots,r_{2m}\in[0, 1]$,
\begin{eqnarray}\label{orbit}
 1_{\{\forall i\neq j, r_i\neq r_j\}} = \sum_{\tau\in\mathfrak{S}_{2m}}
 1_{\{r_{\tau(1)}>\ldots>r_{\tau(2m)}\}} &=&
 \sum_{o_i\in\mathcal{O}}\sum_{\tau\in o_i}
 1_{\{r_{\tau(1)}>\ldots>r_{\tau(2m)}\}}\nonumber \\&\leq&
 \prod_{i=1}^{\frac{(2m)!}{m!}}\prod_{j=1}^{m}1_{\{r_{2\tau_i(j)-1}>r_{2\tau_i(j)}\}.}
\end{eqnarray}

Then for any $\eps>0$, we have
 \begin{eqnarray*}E(F_{\eps}^2)&=&\left|\int_{[0,1] ^{2n}}f_1(x_{1})\ldots f_{2n}(x_{2n})1_{\{\vert
x_{1}-x_{2}\vert <\eps\}}1_{\{\vert x_{1}- x_{3}\vert
<\eps\}}1_{\{\vert x_{n+1}-x_{n+2}\vert <\eps\}}1_{\{\vert x_{n+1}-
x_{n+3}\vert <\eps\}}\right.\\&&\times
\left.E\left(\prod_{j=1}^{n}\theta _{\eps} (x_{j} )\theta _{\eps}
(x_{j+1} ) \right)h_{\eps}(x_{1},\ldots, x_{n})
h_{\eps}(x_{j+1},\ldots, x_{2n})dx_{1}\ldots dx_{2n}\right|\\&\leq&
\int_{[0,1] ^{2n}}|f_1(x_{1})|\ldots |f_{2n}(x_{2n})|1_{\{\vert
x_{1}-x_{2}\vert <\eps\}}1_{\{\vert x_{1}- x_{3}\vert <\eps\}}
\left|E\left(\prod_{j=1}^{2n}\theta _{\eps} (x_{j}
)\right)\right|dx_{1}\ldots dx_{2n}
\\&=& \sum_{o_i\in
\mathcal{O}}\sum_{\tau\in o_i }\int_{[0,1] ^{2n}}|f_1(x_{1})|\ldots
|f_{2n}(x_{2n})|1_{\{\vert x_{1}-x_{2}\vert <\eps\}}1_{\{\vert
x_{1}- x_{3}\vert <\eps\}}\\&&\times
1_{\{x_{\tau(1)}>\ldots>x_{\tau(2n)}\}}Q_{\eps}\left(\sum_{j=1}^{n}(x_{\tau(2j-1)}-x_{\tau(2j)})\right)dx_{1}\ldots
dx_{2n}.
\end{eqnarray*}

Among the addends of the last term there are two possible
situations.

\begin{itemize}
\item On one hand we have terms of the type:
\begin{eqnarray*}
&&\int_{[0,1] ^{4}}1_{\{x_3>x_4\}}
|f_{2\tau_i(k)-1}(x_{1})||f_{2\tau_i(k)}(x_{2})|Q_{\eps}\left(x_1-x_2\right)\\&\times&
1_{\{0<x_1-x_2<\eps\}}1_{\{0<x_1-x_3<\eps\}}
|f_{2\tau_i(k')-1}(x_{3})||f_{2\tau_i(k')}(x_{4})|Q_{\eps}\left(x_3-x_4\right)dx_1dx_2dx_3dx_4\\&\times&
\prod_{j\neq k,k'; =1}^{n}\int_{[0,1] ^{2}}1_{\{x_1>x_2\}}
|f_{2\tau_i(j)-1}(x_{1})||f_{2\tau_i(j)}(x_{2})|Q_{\eps}\left(x_1-x_2\right)dx_1dx_2,
\end{eqnarray*}where $\tau_i(k)>\tau_i(k')+1$.

Notice that, using (\ref{orbit})(as in \cite{BNRT}), we obtain
 \begin{eqnarray*}&&\int_{[0,1] ^{2}}1_{\{x_1>x_2\}}
|f_{2\tau_i(j)-1}(x_{1})||f_{2\tau_i(j)}(x_{2})|Q_{\eps}\left(x_1-x_2\right)dx_1dx_2\\&&\leq
\frac{1}{2} \|f_{2\tau_i(j)-1}\|_{L^2}\|f_{2\tau_i(j)}\|_{L^2}.
\end{eqnarray*}

Moreover, given $h_i\in L^2([0,1])$, $i=1,\dots,4$, we have that
\begin{eqnarray*}
&&\int_{[0,1]^4}|h_1(x_1)||h_2(x_2)||h_3(x_3)||h_4(x_4)|1_{\{0<x_1-x_2<\eps\}}1_{\{0<x_1-x_3<\eps\}}\\&&\times Q_{\eps}\left(x_1-x_2\right)Q_{\eps}\left(x_3-x_4\right)dx_1dx_2dx_3dx_4\\
&\leq& \int_{[0,1]^4}|h_1(x_1)||h_2(x_2)||h_3(x_3)||h_4(x_4)|1_{\{0<x_1-x_2<\eps\}}1_{\{0<x_1-x_3<\eps\}}1_{\{0<x_3-x_4<\eps\}}\\&&\times Q_{\eps}\left(x_1-x_2\right)Q_{\eps}\left(x_3-x_4\right)dx_1dx_2dx_3dx_4\\
&&+ \int_{[0,1]^4}|h_1(x_1)||h_2(x_2)||h_3(x_3)||h_4(x_4)|1_{\{0<x_1-x_2<\eps\}}1_{\{0<x_1-x_3<\eps\}}1_{\{\eps<x_3-x_4\}}\\&&\times Q_{\eps}\left(x_1-x_2\right)Q_{\eps}\left(x_3-x_4\right)dx_1dx_2dx_3dx_4\\
\\&:=&A^1_{\eps}+A^2_{\eps}.
\end{eqnarray*}
The term $A^1_{\eps}$ converges to zero by using the same manner of
the convergence of $I^{\eps}_{2}$ in the proof of Lemma \ref{imp}.
For the term $A^2_{\eps}$ we have that
\begin{eqnarray*}
&&A^2_{\eps}\\&\leq&2\left(\int_{[0,1]^4}h_1^2(x_1)h_3^2(x_3)1_{\{0<x_1-x_3<\eps\}} Q_{\eps}\left(x_1-x_2\right)Q_{\eps}\left(x_3-x_4\right)dx_1dx_2dx_3dx_4\right.\\
&&+\left.\int_{[0,1]^4}h_2^2(x_2)h_4^2(x_4)1_{\{0<x_1-x_3<\eps\}}1_{\{\eps<x_3-x_4\}}
Q_{\eps}\left(x_1-x_2\right)Q_{\eps}\left(x_3-x_4\right)dx_1dx_2dx_3dx_4\right).
\end{eqnarray*}

Integrating with respect to $x_2$ and $x_4$ in the first addend we
obtain the convergence to zero by using the dominated convergence
theorem. Moreover, using the fact that for $y>\eps $,
$Q_{\eps}(y)\leq e^{-2}$, and integrating after with respect to
$x_3$ and $x_1$ we can bound the second addend by
$$\eps\int_{[0,1]^2}h_2^2(x_2)h_4^2(x_4)dx_2dx_4$$
that clearly converge also to zero.

\item We have also terms of the type:
\begin{eqnarray*}
&&\int_{[0,1] ^{6}}1_{\{x_1>x_4\}}
|f_{2\tau_i(k)-1}(x_{1})||f_{2\tau_i(k)}(x_{4})|Q_{\eps}\left(x_1-x_4\right)\\&\times&1_{\{x_2>x_5\}}
|f_{2\tau_i(k')-1}(x_{2})||f_{2\tau_i(k')}(x_{5})|Q_{\eps}\left(x_2-x_5\right)\\&\times&1_{\{x_3>x_6\}}
|f_{2\tau_i(k'')-1}(x_{3})||f_{2\tau_i(k'')}(x_{6})|Q_{\eps}\left(x_3-x_6\right)\\&\times&
1_{\{0<x_1-x_2<\eps\}}1_{\{0<x_1-x_3<\eps\}}dx_1dx_2dx_3dx_4dx_5dx_6\\&\times&
\prod_{j\neq k,k',k''; =1}^{n}\int_{[0,1] ^{2}}1_{\{x_1>x_2\}}
|f_{2\tau_i(j)-1}(x_{1})||f_{2\tau_i(j)}(x_{2})|Q_{\eps}\left(x_1-x_2\right)dx_1dx_2,
\end{eqnarray*}where $\tau_i(k)>\tau_i(k')+1>\tau_i(k'')+2$.

But, using  arguments similar to those presented  in the previous situation it is
not difficult to see that also this type of terms converges to zero.

\end{itemize}

Combining the above convergences we conclude that $E(F_{\eps}^2)$
converges to zero and thus the Lemma \ref{convR} satisfied.

{\bf 2)} Assume now that $ \theta _{\eps}$  are Donsker kernels. For
any $m\geq3$

\begin{eqnarray*}G_{\eps,m}(x_1,\ldots,x_m):=\frac{E(\xi_1^m)}{\varepsilon^m}\sum_{k=1}^{\infty}
1_{[(k-1)\varepsilon^2,k\varepsilon^2)^m}(x_1,\ldots,x_m).
\end{eqnarray*}
 Fix $x_1$  in $[0,1]$. Then for
every $\varepsilon>0$ close to zero, there exist
$k(x_1,\varepsilon)\in\{1,\ldots,\left[\frac{1}{\varepsilon^2}\right]+1\}$
such that $(k(x_1,\varepsilon)-1)\varepsilon^2\leq
x_1<k(x_1,\varepsilon)\varepsilon^2$,  this implies that
$(k(x_1,\varepsilon)-1)\varepsilon^2\rightarrow x_1$ as $\varepsilon
\rightarrow0$. Then we can write
\begin{eqnarray*}J_{\eps,m}&:=&\int_{[0,1]^{m}}f_1(x_1)f_{2}(x_{2})\ldots f_m(x_{m})G_{\eps,m}(x_1,\ldots,x_m)dx_{1}\ldots
dx_{m}\\&=&
E(\xi_1^m)\eps^{m-2}\int_0^1f_1(x_1)\prod_{j=2}^{m}\left[\frac{1}{\varepsilon^2}
\int_{(k(x_1,\varepsilon)-1)\varepsilon^2}^{k(x_1,\varepsilon)\varepsilon^2}f_{j}(x_{j})dx_{j}\right]dx_1.\end{eqnarray*}
 Moreover for each $j=2,\ldots,m$,  the term $\frac{1}{\varepsilon^2}\int_{(k(x_1,\varepsilon)-1)\varepsilon^2}^{k(x_1,\varepsilon)\varepsilon^2}f_{j}(x_{j})dx_{j}$
 converges to $f_j(x_1)$. Combining this with $m\geq3$, we obtain
 that
\begin{eqnarray}\label{convergence of Jm} J_{\eps,m}\underset{\varepsilon
\rightarrow0}{\longrightarrow}0
\end{eqnarray}
On the other hand, if we denote by

\begin{eqnarray*}\bar{G}_{\sigma, \eps}(x_1,\ldots,x_{2n})=\frac{1}{\varepsilon^{2n}}\prod_{j=0}^{n-1}\left(\sum_{k=1}^{\infty}
1_{[(k-1)\varepsilon^2,k\varepsilon^2)^2}(x_{\sigma(2j+1)},x_{\sigma(2j+2)})\right)
\end{eqnarray*}
Fix  $x_{\sigma(2j+1)}\in[0,1]$ for any $j=0,\ldots,n-1$. Then for
every $\varepsilon>0$ close to zero, there exist
$k(x_{\sigma(2j+1)},\varepsilon)\in\{1,\ldots,\left[\frac{1}{\varepsilon^2}\right]+1\}$
such that $$k(x_{\sigma(2j+1)},\varepsilon)\neq
k(x_{\sigma(2j'+1)},\varepsilon)\ \forall j'\neq j \mbox{ and }
(k(x_{\sigma(2j+1)},\varepsilon)-1)\varepsilon^2\leq
x_{\sigma(2j+1)}<k(x_{\sigma(2j+1)},\varepsilon)\varepsilon^2$$ this
implies that
$(k(x_{\sigma(2j+1)},\varepsilon)-1)\varepsilon^2\rightarrow
x_{\sigma(2j+1)}$  as $\varepsilon \rightarrow0$. Then we can write
\begin{eqnarray*} \bar{J}_{\sigma,
\eps}&:=&\int_{[0,1]^{2n}}f_1(x_{1})\ldots
f_{2n}(x_{2n})\bar{G}_{\sigma, \eps}(x_1,\ldots,x_{2n})1_{\{\vert
x_{1}-x_{2}\vert <\eps\}}1_{\{\vert x_{1}- x_{3}\vert
<\eps\}}dx_{1}\ldots
dx_{2n}\\&\leq&\int_{[0,1]^{n}}\prod_{l=0}^{n-1}\left[f_{{\sigma(2j+1)}}(x_{\sigma(2j+1)})\frac{1}{\varepsilon^2}\int_{(k(x_{\sigma(2j+1)},\varepsilon)-1)
\varepsilon^2}^{k(x_{\sigma(2j+1)},\varepsilon)\varepsilon^2}
1_{\{\vert x_{1}-x_{2}\vert <\eps\}}1_{\{\vert x_{1}- x_{3}\vert
<\eps\}}\right.\\&&\times\left.
f_{{\sigma(2j+2)}}(x_{\sigma(2j+2)})dx_{\sigma(2j+2)}\right]
dx_{\sigma(1)}dx_{\sigma(3)}\ldots dx_{\sigma(2n-1)}.\end{eqnarray*}
Moreover, this last term converges to
\begin{eqnarray}\label{convergence of J}&&\int_{[0,1]^{n}}1_{\{
x_{\sigma(2k+1)}=x_{\sigma(2k'+1)};\,\textrm{\footnotesize for some
$k\neq k'$}\}}\nonumber\\&\times&
\prod_{l=0}^{n-1}\left[f_{\sigma(2l+1)}(x_{\sigma(2l+1)})f_{\sigma(2l+2)}(x_{\sigma(2l+1)})\right]dx_{\sigma(1)}dx_{\sigma(3)}\ldots
dx_{\sigma(2n-1)}=0.\end{eqnarray}  From (\ref{convergence of Jm}),
(\ref{convergence of J}) and the fact that the term
$E\left(\prod_{j=1}^{n}\theta _{\eps} (x_{j} ) \right)$  is written
as a sum of terms of type $G_{\eps,m}$ or $\bar{G}_{\sigma, \eps}$
we conclude that
\begin{eqnarray*}E(F_{\eps}^2)\underset{\varepsilon
\rightarrow0}{\longrightarrow}0
\end{eqnarray*} \qed

We can state now our approximation result for multiple fractional integrals when the integrand is a simple function.
\begin{prop}
Let $f$ be a simple function of the form (\ref{step}). Then the
finite dimensional distribution of the process (\ref{ieta}) converge
as $\eps \to 0$ to the finite dimensional distributions of $\left(
I_{n}^{H}(f1_{[0,t]}^{\otimes n})\right) _{t\in [0,1]}$.
\end{prop}
{\bf Proof: }If $f$ is a simple function of the form (\ref{step}) then for every $t\in [0,1]$
\begin{eqnarray*}
I_{n}^{H}(f1_{[0,t]}^{\otimes n} )&=& \sum_{k=1}^{m} \alpha _{k} I_{n}^{H} \left(  1_{(a_{k}^{1}, b_{k}^{1}]}\times \ldots 1_{(a_{k}^{n}, b_{k}^{n}]}1_{[0,t] }^{\otimes n} \right)\\
&=& \sum_{k=1}^{m} \alpha _{k} \big(\prod_{i=1}^nI_1^H(1_{(a_{k}
^{i}, b_{k}^{i}]}1_{[0,t]} ) \\&& +\sum_{l=1}^{[n/2]}(-1) ^{l}
\sum_{\footnotesize\begin{array}{c}j_{1}, \dots, j_{2l}=1;\\ j_{i}\,
distinct\end{array}}^{n}\left( \prod _{u\in \{1,\dots,n\}\setminus
\{j_{1},\dots, j_{2l}\} }
I_{1}^{H}\left( 1_{(a_{k} ^{u}, b_{k}^{u}]}1_{[0,t]} \right) \right) \\
&&\times \langle 1_{(a_{k}^{j_{1}}, b_{k}^{j_{1}}]}1_{[0,t]},
1_{(a_{k}^{j_{2}}, b_{k}^{j_{2}}]}1_{[0,t]} \rangle
_{{\cal{H}}}\cdots\langle 1_{(a_{k}^{j_{2l-1}},
b_{k}^{j_{2l-1}}]}1_{[0,t]}, 1_{(a_{k}^{j_{2l}},
b_{k}^{j_{2l}}]}1_{[0,t]} \rangle _{{\cal{H}}}\big).
\end{eqnarray*}
 The approximation $I_{n_{\eps}}(f)_{t} $ can be expressed as

\begin{eqnarray*}
I_{n_{\eps }}(f)_{t}&=& \int_{[0,1] ^{n}} \left( \Gamma _{H}^{(n)} f 1_{[0,t]}^{n}\right) (x_{1}, \dots, x_{n}) \prod _{j=1}^{n} \theta _{\eps }(x_{j}) \left( \prod _{i,j=1; i\not= j} (1- 1_{\{\vert x_{i}-x_{j}\vert <\eps \}}\right) dx_{1} dx_{n}\\
&=& \int_{[0,1] ^{n}} \left( \Gamma _{H}^{(n)} f 1_{[0,t]}^{n}\right) (x_{1}, \dots, x_{n}) \prod _{j=1}^{n} \theta _{\eps }(x_{j})\\
 &&\times \left(1+\sum_{l=1}^{[n/2]} (-1) ^{l} \left( \!\!\!\!\!\sum_{\footnotesize\begin{array}{c}k_{1},\dots, k_{2l} =1;\\ k_{j}\, distinct\end{array} }^{n} \!\!\!\!\! 1_{\{\vert x_{k_{1}}-x_{k_{2}\vert}<\eps \}} \cdots 1_{\{\vert x_{k_{2l-1}}-x_{k_{2l}\vert}<\eps \}}dx_{1}\cdots dx_{n}\right)\right) + R.
\end{eqnarray*}
The term $R$ above contains terms of the type
\begin{equation*}
 \int_{[0,1] ^{n} }\left( \prod_{j=1}^{n}(\Gamma _{H}^{(1)} 1_{(a_{k}^{j}, b_{k}^{j}]}1_{[0,t]})(x_{j} )\theta_{\eps }(x_{j}) \right) 1_{\{\vert x_{1}-x_{2}\vert <\eps \}} 1_{\{\vert x_{1}-x_{3}\vert <\eps\} } 1_{A}(x_{1},\dots,x_{n})
   \end{equation*}where $A$ is a Borel subset of $[0,1]^{\otimes n}$. It will converge to zero by using Lemma \ref{convR} for $h_{\eps}(x_{1},\dots, x_{n}) =1_{A}(x_{1},\dots, x_{n}).$
The behavior of $I_{n_{\eps}}(f)_{t}$ will be then given by the behavior of
\begin{eqnarray*}
&&\int_{[0,1] ^{n}} \left( \Gamma _{H}^{(n)} f 1_{[0,t]}^{n}\right)
(x_{1}, \dots, x_{n}) \prod _{j=1}^{n} \theta _{\eps
}(x_{j})\\&\times&\left(1+\sum_{l=1}^{[n/2]} (-1) ^{l} \left(
\!\!\!\!\!\sum_{\footnotesize\begin{array}{c}k_{1},\dots, k_{2l}
=1;\\ k_{j}\, distinct\end{array} }^{n} \!\!\!\!\! 1_{\{\vert
x_{k_{1}}-x_{k_{2}\vert}<\eps \}} \cdots 1_{\{\vert
x_{k_{2l-1}}-x_{k_{2l}\vert}<\eps \}}dx_{1}\cdots
dx_{n}\right)\right) .
 \end{eqnarray*}
 First we note that by Lemma 2 the first term in the above expression converges  in the sense of finite dimensional distributions to
\begin{equation*}
\sum_{k} \alpha _{k} I_{1}^{H} (1_{(a_{k}^{1}, b_{k}^{1}]}1_{[0,t]
}) \cdots I_{1}^{H} (1_{(a_{k}^{n}, b_{k}^{n}]}1_{[0,t] })=\sum_{k}
\alpha _{k}(B^H_{b_{k}^{1}\wedge t} -B^H_{a_{k}^{1}\wedge
t})\cdots(B^H_{b_{k}^{n}\wedge t} -B^H_{a_{k}^{n}\wedge t}).
\end{equation*}

 We will show that for every $l=1, \ldots , [\frac{n}{2}] $ and for every $j_{1}, \ldots , j_{2l} =1,\dots, n$ distinct the sequence
\begin{equation*}
\int_{[0,1] ^{n}}dx_{1}\cdots dx_{n} \Gamma _{H} ^{(n)} (f1_{[0,t]
}^{\otimes n} )(x_{1}, \dots, x_{n}) \left( \prod_{j=1} ^{n} \theta
_{\eps }(x_{j})\right) 1_{\{\vert x_{j_{1} } - x_{j_{2}}\vert <\eps
\}} \cdots 1_{\{\vert x_{j_{2l-1} } - x_{j_{2l}}\vert <\eps \}}
\end{equation*}
converges in the sense of finite dimensional distributions to the stochastic process
\begin{eqnarray*}
&&\sum_k \alpha_k\left( \prod_{u\in \{1,\dots, n\}\setminus \{
j_{1}, \dots, j_{2l}\} } I_{1}^{H} \left( 1_{(a_{k}^{u}, b_{k} ^{u}
]} 1_{[0,t]} \right) \right)\\&\times&
 \langle 1_{(a_{k}^{j_{1}}, b_{k}^{j_{1}}]}1_{[0,t]}, 1_{(a_{k}^{j_{2}}, b_{k}^{j_{2}}]}1_{[0,t]} \rangle _{{\cal{H}}}\cdots\langle 1_{(a_{k}^{j_{2l-1}}, b_{k}^{j_{2l-1}}]}1_{[0,t]}, 1_{(a_{k}^{j_{2l}}, b_{k}^{j_{2l}}]}1_{[0,t]} \rangle _{{\cal{H}}}.
\end{eqnarray*}
Indeed, since
\begin{equation*}
\Gamma _{H} ^{(n)} (f1_{[0,t] }^{\otimes n} )(x_{1}, \dots, x_{n}) =
\sum_{k=1}^{m} \alpha _{k} (\Gamma_{H}^{(1)} 1_{(a_{k}^{1}, b_{k}
^{1}]})(x_{1})\cdots (\Gamma_{H}^{(1)} 1_{(a_{k}^{n}, b_{k}
^{n}]})(x_{n})
\end{equation*}
we can write, for every $j_{1},\dots, j_{2l} =1,\dots, n$ distinct
\begin{eqnarray*}
&&\int_{[0,1] ^{n}}dx_{1}\cdots dx_{n} \Gamma _{H} ^{(n)} (f1_{[0,t] }^{\otimes n} )(x_{1}, \dots, x_{n}) \left( \prod_{j=1} ^{n} \theta _{\eps }(x_{j})\right) 1_{\{\vert x_{j_{1} } - x_{j_{2}}\vert <\eps \}} \cdots 1_{\{\vert x_{j_{2l-1} } - x_{j_{2l}}\vert <\eps\} }\\
&=& \sum_{k=1}^{m} \alpha _{k} \left(   \prod_{u\in \{1,\dots, n\}\setminus \{j_{1},\dots, j_{2l}\} }\int_{[0,1] }dx_{u} (\Gamma _{H} ^{(1)} 1_{(a_{k}^{u}, b_{k}^{u}]}1_{[0,t] } )(x_{u}) \theta _{\eps}(x_{u})\right)\\
&&\times \int_{[0,1]^{2}} (\Gamma _{H}^{(1)} 1_{(a_{k}^{j_{1}}, b_{k} ^{j_{1}}]} 1_{[0,t]} )(x_{j_{1}})(\Gamma _{H}^{(1)} 1_{(a_{k}^{j_{2}}, b_{k} ^{j_{2}}]} 1_{[0,t]} )(x_{j_{2}})\theta _{\eps}(x_{j_{1}}) \theta _{\eps}(x_{j_{2}}) 1_{\{\vert x_{j_{1}}-x_{j_{2}} \vert <\eps\}}dx_{j_{1}}dx_{j_{2}}\\
&&\ldots \times \int_{[0,1]^{2}} (\Gamma _{H}^{(1)} 1_{(a_{k}^{j_{2l-1}}, b_{k} ^{j_{2l-1}}]} 1_{[0,t]} )(x_{j_{2l-1}})(\Gamma _{H}^{(1)} 1_{(a_{k}^{j_{2l}}, b_{k} ^{j_{2l}}]} 1_{[0,t]} )(x_{j_{2l}})\\
&&\hskip2cm \times \theta _{\eps}(x_{j_{2l-1}}) \theta _{\eps}(x_{j_{2l}}) 1_{\{\vert x_{j_{2l-1}}-x_{j_{2l}} \vert <\eps\}}dx_{j_{2l-1}}dx_{j_{2l}}\\
&=& \sum_{k=1}^{m} \alpha _{k} \left(   \prod_{u\in \{1,\dots,
n\}\setminus \{j_{1},\dots, j_{2l}\} }(\eta _{\eps} (b_{k}^{u}\wedge
t )- \eta _{\eps} (a_{k}^{u} \wedge t)\right)\\&&\times \left(
Y^{\eps}_{b_{k}^{j_{1}}\wedge t, b_{k}^{j_{2}}\wedge t}
-Y^{\eps}_{a_{k}^{j_{1}}\wedge t, b_{k}^{j_{2}}\wedge
t}-Y^{\eps}_{b_{k}^{j_{1}}\wedge t, a_{k}^{j_{2}}\wedge
t}+Y^{\eps}_{a_{k}^{j_{1}}\wedge t, a_{k}^{j_{2}}\wedge t}\right)
 \end{eqnarray*}
where, for $v=1,\dots, l$ we denoted by
\begin{eqnarray*}
&& Y^{\eps}_{s,t}  = \int_{[0,1]^{2}} (\Gamma _{H}^{(1)}  1_{[0,t]}
)(x_{1})(\Gamma _{H}^{(1)} 1_{[0,s]} )(x_{2}) \theta _{\eps}(x_{1})
\theta _{\eps}(x_{2}) 1_{\{\vert x_{1}-x_{2} \vert
<\eps\}}dx_{1}dx_{2}.
\end{eqnarray*}
The conclusion follows by using Lemma \ref{lema6} and the results
obtained for the case $n=2$. \qed

\vskip0.3cm

We state now our main result.

\begin{theorem}
Let $f$ be a  function in the space $\left| {\cal{H}}\right| ^{\otimes n}$. Then the finite dimensional distribution of the process (\ref{ieta}) converge as $\eps \to 0$ to the finite dimensional distributions of $\left( I_{n}^{H}(f1_{[0,t]}^{\otimes n})\right) _{t\in [0,1]}$.
\end{theorem}
{\bf Proof: } It is a consequence of Lemma 2.1 and Theorem 2.3 in
\cite{BJT2}, of the isometry of multiple integrals and of the fact
that the simple functions are dense in  $\left| {\cal{H}}\right| ^{\otimes n}$
since for every $t\in [0,1] $ it holds (see Section 2.2 in
\cite{BJT2})
\begin{equation*}
E\left| I_{n_{\eps}}(f)_{t} \right| \leq c\Vert \Gamma _{H}^{(n)}f 1_{[0,t]}^{\otimes n} \Vert _{L^{2}([0,1] ^{\otimes n})} =c  \Vert f 1_{[0,t]}^{\otimes n} \Vert _{{\cal{H}}^{\otimes n}}.
\end{equation*}
\qed

{\bf Acknowledgement: } Part of this paper has been written while
the third author was visiting the research center ``Centre de
Recerca Matem\`atica" in Barcelona in February 2009. He warmly
acknowledges support and hospitality.

\end{document}